\documentclass[aap,MSNbibl,seceqn,dvips]{arximspdf}

\doi{10.1214/13-AAP956} %kopijuoti is PTS
\volume{24}
\issue{4}
\pubyear{2014}
\firstpage{1554}
\lastpage{1584}

\makeatletter

\def\cal{\mathcal}
\newtheorem{theorem}{Theorem}[section]
\newtheorem{lemma}[theorem]{Lemma}
\newtheorem{lemmas}{Lemma}[section]
\newtheorem{corollary}[theorem]{Corollary}
\newtheorem{proposition}[theorem]{Proposition}
\newproclaim{definition}[theorem]{Definition}
\newproclaim{remark}[theorem]{Remark}
\newproclaim{example}[theorem]{Example}

\def\bfR{\mathbb{R}}
\def\bfS{\mathcal{S}}
\newcommand{\esssup}{\mathop{\operatorname{ess\,sup}}}
\makeatother

\begin{document}
\begin{frontmatter}

\title{Monotonicity of the value function
for a two-dimensional optimal stopping problem}

\runtitle{Monotonicity value function}

\begin{aug}
\author{\fnms{Sigurd} \snm{Assing}\corref{}\ead[label=e1]{s.assing@warwick.ac.uk}},
\author{\fnms{Saul} \snm{Jacka}\ead[label=e2]{s.d.jacka@warwick.ac.uk}}
\and
\author{\fnms{Adriana} \snm{Ocejo}\thanksref{t1}\ead[label=e3]{a.ocejo-monge@warwick.ac.uk}}
\thankstext{t1}{Supported by CONACYT Ph.D. scholarship \#309247.}
\runauthor{S. Assing, S. Jacka and A. Ocejo}
\affiliation{University of Warwick}
\address{Department of Statistics\\
University of Warwick\\
Coventry CV4 7AL\\
United Kingdom\\
\printead{e1}\\
\phantom{E-mail:\ }\printead*{e2}\\
\phantom{E-mail:\ }\printead*{e3}}
\end{aug}

% HISTORY:
\received{\smonth{8} \syear{2012}}
\revised{\smonth{7} \syear{2013}}

% ABSTRACT

\begin{abstract}
We consider a pair $(X,Y)$ of stochastic processes satisfying the equation
$dX=a(X)Y\,dB$
driven by a Brownian motion and study the monotonicity and
continuity in $y$ of the value function
$v(x,y)=\sup_{\tau} E_{x,y} [ e^{- q\tau}g(X_\tau)]$,
where the supremum is taken over stopping times with respect to
the filtration generated by $(X,Y)$.
Our results can successfully be applied to pricing
American options where $X$ is the discounted price of an asset
while $Y$ is given by a stochastic volatility model such
as those proposed by Heston or Hull and White.
The main method of proof is based on time-change and coupling.
\end{abstract}

% KEYWORDS
\begin{keyword}[class=AMS]
\kwd[Primary ]{60G40}
\kwd[; secondary ]{91G20}
\end{keyword}
\begin{keyword}
\kwd{Optimal stopping}
\kwd{time-change}
\kwd{coupling}
\kwd{stochastic volatility model}
\kwd{American option}
\end{keyword}

\end{frontmatter}

%s1 #&#
\section{Introduction}
Consider a two-dimensional strong Markov process
$(X,\break Y)=(X_t,Y_t, t\geq0)$
with state space $\bfR\times\bfS, \bfS\subseteq(0,\infty)$,
given on a family of probability spaces
$(\Omega, \mathcal{F}, P_{x,y},(x,y)\in\bfR\times\bfS)$
which satisfies the stochastic differential equation
%
%e1.1 #&#
\begin{equation}
\label{eq:dynX} dX = a(X) Y\,dB,
\end{equation}
where $B=(B_t)_{t\geq0}$ is a standard Brownian motion, and
$a\dvtx \bfR\to\bfR$ is a measurable function.

Processes of this type are common in mathematical finance, and in
this context, $X$ would be the discounted price of an asset
while $Y$ is a process giving the so-called stochastic volatility.

We shall refer to this application in the examples, as it was our motivation
in the beginning. However,
the methods used are of a broader nature and can be applied in a wider context.

This paper mainly deals with the regularity of the value function
%
%e1.2 #&#
\begin{equation}
\label{osp} v(x,y)=\sup_{0\leq\tau\leq T} E_{x,y} \bigl[
e^{- q\tau}g(X_\tau)\bigr],\qquad (x,y)\in\bfR\times\bfS,
\end{equation}
with respect to the optimal stopping problem given by $(X,Y)$, a
discount rate
\mbox{$q>0$}, a time horizon $T\in[0,\infty]$ and a measurable
gain function $g\dvtx \bfR\to\bfR$. But for financial applications
(see Section~\ref{option pricing}),
a slightly modified value function of type
{\renewcommand{\theequation}{1.2$'$}
\begin{equation}\label{eq1.2prime}
v(x,y)=\sup_{0\leq\tau\leq T} E_{x,y} \bigl[ e^{- r\tau}g
\bigl(e^{r\tau
}X_\tau\bigr)\bigr] %\leqno(1.2^\prime)
\end{equation}}
\hspace*{-2pt}is also considered where $r$ stands for the instantaneous interest rate.

The supremum in (\ref{osp}) and (\ref{eq1.2prime}) is taken over all
\emph{finite} stopping times with respect to
the filtration generated by the pair of processes $(X,Y)$.

%In the case $T=\infty$ we could define $e^{-qT}g(X_T)=0$ but this

%would still not ensure that there is a stopping time taking the

%supremum.

To ensure the well-posedness of this problem, we assume
the integrability condition (recall that $T$ may be infinite)
%
%e1.3 #&#
\setcounter{equation}{2}
\begin{equation}
\label{eq:integrabilityCond} E_{x,y} \Bigl[ \sup_{0\leq t \leq T}e^{-qt}
\bigl|g(X_t)\bigr| I(t<\infty) \Bigr]< \infty \qquad\mbox{for all } (x,y)\in\bfR\times
\bfS,
\end{equation}
which is a common assumption in the context of optimal stopping problems.

Note that this condition is satisfied if $g$ is bounded. For more
general functions, verifying this condition can be fairly difficult,
and its validity may depend on the particular choice of the dynamics for
$(X,Y)$.

%E_{x,y} \left[ \phantom{\rule}{0pt}{11pt}

%e^{- q\tau}|g(X_\tau)|\right] < \infty

% \mbox{for all} (x,y)\in\bfR\times\bfS

%for all stopping times $\tau$ with respect to

%the filtration generated by $(X,Y)$. Remark that,

%under this condition, the value function might become infinitely large

%but such a singularity would not affect our comparison method.

%Note that (\ref{eq:integrabilityCond}) is trivially satisfied if $g$
%is bounded.

Our main focus is on proving
the monotonicity %and continuity of
of $v(x,y)$ with respect to $y\in\bfS$,
and we are able to verify this property in the case
of the following two classes of strong Markov processes
under not too restrictive conditions
(see Theorems~\ref{main one} and \ref{main two}):
\begin{itemize}\label{classes}
\item
\textit{Regime-switching}: $Y$ is a skip-free
continuous-time Markov chain (see page~\pageref{skippy}) which
is independent of the Brownian motion $B$ driving equation
(\ref{eq:dynX}).
\item
\textit{Diffusion}: $Y$ solves a stochastic differential equation
of the type
%
%e1.4 #&#
\begin{equation}
\label{vol_W} dY = \eta(Y)\,dB^Y + \theta(Y)\,dt,
\end{equation}
where $B^Y=(B^Y_t)_{t\geq0}$ is a standard Brownian motion
such that the quadratic covariation satisfies
$\langle B, B^Y \rangle_t=\delta t, t\ge0$, for some real
parameter $\delta\in[-1,1]$ and $\eta,\theta\dvtx \bfR\to\bfR$ are
measurable functions.
\end{itemize}
Note that, in the second class, the joint distribution of $X$
and $Y$ is uniquely determined if the system of equations
(\ref{eq:dynX}), (\ref{vol_W}) admits a weakly unique solution, and
the process $Y$ does not have
to be independent of the driving Brownian motion $B$,
whereas, in the case of the first class, the process $Y$ is not
given by an equation, and the assumed independence of $Y$ and $B$ is
a natural way of linking $X$ and $Y$ if there is too little
information about the structure of the pair $(X,Y)$.

Our technique is based on time-change and coupling.
Equation (\ref{eq:dynX}) goes back to a volatility model used by
Hobson in \cite{Hobson} who also applies time-change and coupling but
for comparing prices of European options.
As far as we know,
our paper is the first paper dealing with the extra difficulty of
applying this technique in the context of optimal stopping. It should
be mentioned that Ekstr\"om \cite{Ekstrom}, Theorem~4.2, can compare
prices of American options if $Y\equiv1$ in equation (\ref{eq:dynX}) and
$a$ also depends on time. Nevertheless, it seems to be that his method
cannot be applied in the case of nontrivial processes $Y$.

We provide some examples to illustrate the results.
In the case of regime-switching,
we look at the pricing
of perpetual American put options which, for $a(x)=x$,
was studied by Guo and Zhang \cite{G-Z}
for a two-state Markov chain
and by Jobert and Rogers \cite{J-R} for a finite-state Markov chain.
While the former, since the situation is much easier,
gave a closed-form expression for the price,
the latter could only provide a numerical algorithm to approximate
the value function which gives the price of the contract.
It turns out that the algorithm in the case of a chain with many states
can be very time-intensive if
the unknown \textit{thresholds}
which characterize the optimal stopping rule
are not known to be in a specific order when labeled by the different
volatility states \emph{before} the algorithm starts.
However,
based on our result that the value function $v(x,y)$ is monotone in $y$,
we are now able to give conditions under which these \textit{thresholds}
must be in a monotone order.

Ultimately, in the case where $Y$ is a diffusion,
we verify the continuity and monotonicity of the value function
$v(x,y)$ with respect to $y\in\bfS=(0,\infty)$ for two important volatility
models, the Heston \cite{Heston} and the Hull and White \cite{HW} model.
Note that, using entirely different methods, differentiability and monotonicity
in the volatility parameter of European option prices under the Hull
and White model
were studied in \cite{BR1996,RT1997}. The authors of \cite{RT1997}
also showed a connection between the monotonicity in the volatility
parameter and
the ability of an option to complete the market. Another motivation to study
the monotonicity of the value function in the volatility parameter $y$
is that the
numerical solution of the corresponding free-boundary problem becomes a
lot easier
if we know that the continuation region is monotonic
in $y$ and if we know that the corresponding free-boundary is continuous.
Moreover, we will show, in a sequel, under the assumption of continuity,
how to solve a game-theoretic version of the American put problem
corresponding to model
uncertainty for the stochastic volatility.

The structure of this paper is as follows.
In Section~\ref{regswi} the monotonicity of the value function
$v(x,y)$ with respect to
$y\in\bfS=\{y_i\dvtx i=1,2,\ldots,m\}\subseteq(0,\infty)$
is shown in the case of
regime-switching, and the main method is established.
In Section~\ref{diffu} the main method is adapted to the case of a
system of stochastic differential equations
(\ref{eq:dynX}), (\ref{vol_W}) which is the diffusion case,
while in Section~\ref{conti} we use monotonicity to show the continuity
of the value function $v(x,y)$ with respect to $y\in\bfS=(0,\infty)$
in the diffusion case.
In Section~\ref{option pricing} we reformulate our results in the
context of option pricing.
Then all our examples are discussed in
detail in Section~\ref{example section} and, in the \hyperref[app]{Appendix},
we prove auxiliary results and some of the corollaries.

Finally,
it should be mentioned that all our results and proofs would not
change in principle
if the state space of $(X,Y)$ is $\bfR\times\bfS$ with
$\bfS\subseteq(-\infty,0)$ instead of $\bfS\subseteq(0,\infty)$.
The only change in this case [see Corollary \ref{version for negative}(ii)]
would be to \textit{order}: \textit{increasing} becomes \textit{decreasing}.
However, as pointed out in the proof of
Corollary~\ref{version for negative}(ii),
our method cannot be applied to show the monotonicity of
$v(x,y)$ in $y\in\bfS$ if $\bfS$ contains a neighborhood of zero.
We do not know either how to generalize our method to the
nonmartingale case.

%s2 #&#
\section{The regime-switching case}\label{regswi}

Suppose %that
$({X},{Y})=({X}_t,{Y}_t$, $t\geq0)$ is a strong Markov process
given on a family of probability spaces
$({\Omega}, {\mathcal{F}}$, $P_{x,y},(x,y)\in\bfR\times\bfS)$
which satisfies the following conditions:
\begin{longlist}[(C1)]\label{a and b}
\item[(C1)]
The process
$(X,Y)$ is adapted with respect to a filtration ${\mathcal{F}}_t, t\ge
0$, of sub-$\sigma$-algebras of ${\mathcal{F}}$ and,
for every $(x,y)\in\bfR\times\bfS$, there is an
${\mathcal{F}}_t$ Brownian motion ${B}$
on $({\Omega}, {\mathcal{F}}, P_{x,y})$
independent of ${Y}$ such that
\[
{X}_t = {X}_0+ \int_0^t
a({X}_s){Y}_s \,d{B}_s,\qquad t\ge0,
P_{x,y}\mbox{-a.s.};
\]
\item[(C2)]
The process $Y$ is a continuous-time Markov chain
on the finite state space
$\bfS=\{y_i\dvtx i=1,2,\ldots,m\} \subset(0,\infty)$
with $Q$-matrix $(q[y_i,y_j])$.
\end{longlist}
%
%re2.1 #&#
\begin{remark}\label{G strong markov}
(i) Because of the condition $\min\{y_1,\ldots,y_m\}>0$ we have that
$P_{x,y} ( \lim_{t\uparrow\infty}$ $\int_0^t {Y}^2_s \,ds =
\infty)=1$ for all $(x,y)\in\bfR\times\bfS$.

\begin{longlist}[(iii)]
\item[(ii)] From the above assumptions it immediately follows that,
for every initial condition $x\in\bfR$, there exists
a weak solution to the stochastic differential equation
$dG = a(G)\,dW$ driven by a Brownian motion $W$.
To see this fix $(x,y)\in\bfR\times\bfS$, and write
\[
{X}_t = x+ \int_0^t
a({X}_s) \,dM_s,\qquad t\ge0, P_{x,y}\mbox{-a.s.},
\]
where $M_s =\int_0^s Y_u \,dB_u$ is well defined since
$\int_0^s Y_u^2\,du < \infty, P_{x,y}$-a.s., for all $s\ge0$.
But time-changing $X$ by the inverse of $\langle M\rangle$,
which exists by (i) above, yields
\[
{G}_t = x+ \int_0^t
a({G}_s) \,dW_s, \qquad t\ge0, P_{x,y}\mbox{-a.s.},
\]
where $G=X\circ\langle M\rangle^{-1}$ is
${\mathcal{F}}_{\langle M\rangle^{-1}_t}$-adapted,
and $W=M\circ\langle M\rangle^{-1}$ is an
${\mathcal{F}}_{\langle M\rangle^{-1}_t}$ Brownian motion by
the Dambis--Dubins--Schwarz theorem; see \cite{Rev-Yor},\break  Theorem~V.1.6.
The equation does indeed hold for all $t\ge0$ since\break
$P_{x,y} ( \lim_{t\uparrow\infty} \int_0^t {Y}^2_s \,ds = \infty)=1$.

\item[(iii)] Because $\langle M\rangle^{-1}
=\int_0^\cdot{Y}_{\langle M\rangle^{-1}_s}^{-2} \,ds$,
an easy calculation shows that the process
$Y\circ\langle M\rangle^{-1}$ is a continuous-time Markov chain
with $Q$-matrix $(y_i^{-2}q[y_i,\break y_j]),  y_i,y_j\in\bfS$.
%which takes the same values $\bfS=\{y_i:i=1,2,\ldots\} \subset(0,
\end{longlist}
\end{remark}

We can now formulate the condition on the coefficient $a$ needed for
our method.
\begin{longlist}[(C3)]
\item[(C3)]\label{new C3}
Let $a\dvtx \bfR\to\bfR$ be measurable functions such that the stochastic differential equation
$dG = a(G)\,dW$ driven by a Brownian motion $W$
has a weakly unique strong Markov solution with state space $\bfR$.
\end{longlist}

The law of the strong Markov process given by {(C3)} is entirely
determined
by its semigroup of transition kernels.
Multiplying these transition kernels and
the transition kernels of a continuous-time Markov chain on $\bfS
\times
\bfS$
both marginals of which are determined by the $Q$-matrix
$(y_i^{-2}q[y_i,y_j]), y_i,y_j\in\bfS$,
results in a semigroup of transition kernels of a strong Markov process
$(G,Z,Z^\prime)$ with $G$ being independent of $(Z,Z^\prime)$.
Now choose a complete probability space
$(\tilde{\Omega}, \tilde{\mathcal{F}},\tilde{P})$
such that $(G,Z,Z^\prime)$ starts from fixed
$(x,y,y^\prime)\in\bfR\times\bfS\times\bfS$.
Let ${\cal F}^{G,Z,Z^\prime}_{t}$ denote the augmentation of the filtration
$\sigma(\{G_s,Z_s,Z^\prime_s\dvtx\break  s\le t\}), t\ge0$,
and assume that
$(G,Z),(G,Z^\prime)$
are strong Markov processes
with respect to ${\cal F}^{G,Z,Z^\prime}_{t}$---an example will be
given in the proof of Theorem~\ref{main one}.

Moreover, by the
martingale problem associated with the strong Markov process $G$,
$G_t-x$
is a continuous local ${\cal F}^{G,Z,Z^\prime}_{t}$-martingale
with quadratic variation $\int_0^t a(G_s)^2\,ds, t\ge0$.
Thus, by a well-known result going back to Doob
(see \cite{Ik-Wa}, Theorem~II 7.1$'$, e.g.),
there is a Brownian motion $W$ such that
%
%e2.1 #&#
\begin{equation}
\label{G-W-equ} G_t-x =\int_0^t
a(G_s)\,dW_s,\qquad t\ge0, \tilde{P}\mbox{-a.s}.
\end{equation}
The construction of $W$ on
$(\tilde{\Omega},\tilde{\mathcal{F}},\tilde{P})$
(or on a canonical enlargement of it\setcounter{footnote}{1}\footnote{Our convention
is to use $(\tilde{\Omega},\tilde{\mathcal{F}},\tilde{P})$
for the enlarged space, too.})
as given in the proof of Theorem II 7.1$'$ in \cite{Ik-Wa} shows that
the pair $(G,W)$ is also independent of $(Z,Z^\prime)$.
But note that $W$ might only be a Brownian motion with respect to a
filtration $\tilde{\cal F}_t$ larger than
${\cal F}^{G,Z,Z^\prime}_{t}, t\ge0$, so that the stochastic integral
in (\ref{G-W-equ}) can only be understood with respect to the larger
filtration.
%
%co2.2 #&#
\begin{corollary}\label{basic}
For given $(x,y,y^\prime)\in\bfR\times\bfS\times\bfS$,
there is a complete probability space
$(\tilde{\Omega}, \tilde{\mathcal{F}},\tilde{P})$
equipped with two filtrations
${\cal F}^{G,Z,Z^\prime}_{t}\subseteq\tilde{\cal F}_t, t\ge0$,
which is big enough to carry four basic processes
$G,W,{Z},{Z}^\prime$ such that:
$(G,W)$ is a weak $\tilde{\cal F}_t$-adapted
solution of $dG = a(G)\,dW$ starting from $x$
independent of $(Z,Z^\prime)$,
the processes ${Z}$ and ${Z}^\prime$ are Markov chains
with $Q$-matrices $(y_i^{-2}q[y_i,y_j]), y_i,y_j\in\bfS$,
starting from $y$ and $y^\prime$, respectively,
and $(G,Z),(G,Z^\prime)$ are strong Markov processes
with respect to ${\cal F}^{G,Z,Z^\prime}_{t}, t\ge0$.
\end{corollary}

The goal of this section is to show that,
under some not too restrictive conditions,
for fixed $x\in\bfR$ and $y,y^\prime\in\bfS$,
%
%e2.2 #&#
\begin{equation}
\label{want to compare} \mbox{if } y\le y^\prime\mbox{ then } v(x,y)\le v
\bigl(x,y^\prime\bigr),
\end{equation}
where the value function $v$ is given by (\ref{osp}).

Choosing $x$ and $y\le y^\prime$,
we will construct two processes
$(\tilde{X},\tilde{Y})$ and $(\tilde{X}^\prime,\tilde{Y}^\prime)$ on
$(\tilde{\Omega}, \tilde{\mathcal{F}},\tilde{P})$ such that
$(\tilde{X},\tilde{Y})$ has the same law as $({X},{Y})$ under $P_{x,y}$,
and $(\tilde{X}^\prime,\tilde{Y}^\prime)$ has the same law as
$({X},{Y})$ under $P_{x,y^\prime}$.
As a consequence we obtain that
%
%e2.3 #&#
\begin{eqnarray}
\label{one prob space} v(x,y)&=& \sup_{0\leq\tilde{\tau}\leq T} \tilde{E} \bigl[
e^{- q\tilde{\tau}}g(\tilde{X}_{\tilde{\tau}})\bigr],
\nonumber
\\[-8pt]
\\[-8pt]
\nonumber% \mbox{and}
v\bigl(x,y^\prime\bigr)&= &\sup_{0\leq\tilde{\tau}^\prime\leq T} \tilde{E} \bigl[
e^{-q\tilde{\tau}^\prime}g\bigl(\tilde{X}^\prime_{\tilde
{\tau
}^\prime}\bigr)\bigr],
\end{eqnarray}
where $\tilde{\tau}$ and $\tilde{\tau}^\prime$
are finite stopping times with respect to the
filtrations generated by $(\tilde{X},\tilde{Y})$ and $(\tilde
{X}^\prime
,\tilde{Y}^\prime)$, respectively.

To see this note that each stopping time ${\tau}$
with respect to the filtration generated by $({X},{Y})$
can easily be associated with a stopping time $\tilde{\tau}$
with respect to the filtration generated by $(\tilde{X},\tilde{Y})$
such that
\[
E_{x,y} \bigl[ e^{- q\tau}g(X_\tau)\bigr] = \tilde{E}
\bigl[ e^{- q\tilde{\tau}}g(\tilde{X}_{\tilde{\tau}})\bigr]
\]
and vice versa proving the first equality in (\ref{one prob space}).
The second equality follows of course by the same argument.

Hence, we can now work on only ONE probability space.
This is an important part of our
method for proving (\ref{want to compare}) which is based
on time-change and coupling and which is demonstrated below.

Let $G,W,{Z},{Z}^\prime$ be given on
$(\tilde{\Omega}, \tilde{\mathcal{F}},\tilde{P})$ as described in
Corollary \ref{basic}, and define
\[
\Gamma_t=\int_0^t
{Z}_s^{-2} \,ds,\qquad t\geq0.
\]
This process $\Gamma=(\Gamma_t)_{t\ge0}$ is of course continuous but
also strictly increasing since $Z$ only takes nonzero values.
Moreover, condition (C2) on page \pageref{a and b} implies that
%
%e2.4 #&#
\begin{equation}
\label{inftyprop} \Gamma_t<\infty,  t\geq0, \mbox{a.s.}\quad\mbox{and}\quad \lim
_{t\uparrow\infty} \Gamma_t = \infty \mbox{ a.s.},
\end{equation}
since, by Remark \ref{G strong markov}(iii), $Z$ has the same law as
$Y\circ\langle M\rangle^{-1}$ under $P_{x,y}$ with
$\int_0^\cdot{Y}_{\langle M\rangle^{-1}_s}^{-2} \,ds$
being the inverse of $\int_0^\cdot Y_s^2\,ds$.
Thus $A=\Gamma^{-1}$
is also a continuous and strictly increasing process satisfying
%
%e2.5 #&#
\begin{equation}
\label{inftyprop for A} A_t<\infty,  t\geq0, \mbox{ a.s.}\quad \mbox{and}\quad \lim
_{t\uparrow\infty} A_t = \infty \mbox{ a.s.}
\end{equation}
As a consequence, the two technical properties:
\begin{longlist}[(P2)]\label{prop_tech}
\item[(P1)]
$\Gamma_{A_t}=A_{\Gamma_t}=t$ for all $t\ge0$ a.s. and
\item[(P2)]
$s<\Gamma_t$ if and only if $A_s<t$
for all $0\leq s,t<\infty$ a.s.
\end{longlist}
must hold.

Of course, $\Gamma$ is adapted to both filtrations
${\cal F}^{G,Z,Z^\prime}_{t}$ and $\tilde{\cal F}_t, t\ge0$.
However, $A=\Gamma^{-1}$ is considered an
${\cal F}^{G,Z,Z^\prime}_{t}$ time change
in the following lemma.
We denote by $\mathcal{M}$ and $\mathcal{T}$ the families of stopping
times with respect to the filtrations
$({\cal F}^{G,Z,Z^\prime}_{t})_{t\geq0}$ and
$({\cal F}^{G,Z,Z^\prime}_{A_t})_{t\geq0}$, respectively.
%
%le2.3 #&#
\begin{lemma} \label{lem_equiv_filt}
If $\rho\in\mathcal{M}$ then $\Gamma_\rho\in\mathcal{T}$,
and
if $\tau\in\mathcal{T}$ then $A_\tau\in\mathcal{M}$.
\end{lemma}
A similar lemma can be found in
{\cite{sharpe}}.
Since the
above lemma is going to be used to reformulate the original optimal
stopping problem (\ref{osp}) in both the case where $Y$ is a Markov
chain and the case where $Y$ is a diffusion, its proof is given in the
\hyperref[app]{Appendix} for
completeness.

The reformulation of (\ref{osp}) is based on the existence of a suitable
solution to~(\ref{eq:dynX}) which is constructed next.

Since $Z$ is $\tilde{\mathcal{F}}_t$-adapted,
one can rewrite (\ref{G-W-equ}) to get
\[
G_t = x+\int_0^t
a(G_s)Z_s\,d\tilde{M}_s, t\ge0,\mbox{a.s.}
\quad\mbox{where}\quad \tilde{M}_s =\int_0^s
\frac{dW_u}{Z_u}, s\ge0.
\]
Observe that the stochastic integral defining $\tilde{M}$ exists
by (\ref{inftyprop}).
Time changing the above equation by $A$ yields
\[
\tilde{X}_t = x+\int_0^t a(
\tilde{X}_s) \tilde{Y}_{_s} \,d\tilde{B}_s,\qquad t
\ge0, \mbox{a.s.}
\]
for $\tilde{X}=G\circ A, \tilde{Y}=Z\circ A,
\tilde{B}=\tilde{M}\circ A$.
Of course, $(\tilde{X},\tilde{Y})$ is $\tilde{\mathcal{F}}_{A_t}$-adapted,
and $\tilde{B}$ is an
$\tilde{\mathcal{F}}_{A_t}$ Brownian motion by
Dambis--Dubins--Schwarz' theorem \cite{Rev-Yor}, Theorem~V.1.6.
Thus $(\tilde{X},\tilde{Y})$ gives a weak solution to (\ref{eq:dynX})
starting from $(x,y)$.
Moreover, $\tilde{B}$ and $\tilde{Y}$ are independent since
$W$ and $Z$ are independent.
The proof of this is contained in the \hyperref[app]{Appendix}; see Lemma A.1 on page
\pageref{A1}.
%
%pr2.4 #&#
\begin{proposition}\label{equiv stop probl}
Let $G,\tilde{X},\tilde{Y}$ be the processes on
$(\tilde{\Omega}, \tilde{\mathcal{F}}, \tilde{P})$
introduced above and starting from $G_0=\tilde{X}_0=x$ and $\tilde{Y}_0=y$.
If the stochastic differential equation
\[
d{X} = a({X}){Y}\,d{B}, \qquad\mbox{$({X},{Y})$ unknown},
\]
driven by a Brownian motion $B$,
where $Y$ is required to be a continuous-time Markov chain
independent of $B$ with $Q$-matrix $(q[y_i,y_j]), y_i,y_j\in\bfS$,
admits a weakly unique solution then,
for any $T\in[0,\infty]$,
\[
v(x,y)=\sup_{\tau\in\mathcal{T}_T } \tilde{E} \bigl[e^{-q \tau} g(
\tilde{X}_\tau)\bigr]=\sup_{\rho
\in\mathcal{M}_T} \tilde{E}
\bigl[e^{-q\Gamma_\rho} g(G_{\rho})\bigr],
\]
where
\[
\mathcal{T}_T=\{\tau\in\mathcal{T} \dvtx 0\leq\tau\leq T\}\quad
\mbox{and}\quad \mathcal{M}_T=\{\rho\in\mathcal{M} \dvtx 0\leq\rho\leq
A_T\}.
\]
Here,
$\mathcal{T}$ and $\mathcal{M}$
denote the families of \emph{finite} stopping times with respect
to the
filtrations
$({\cal F}^{G,Z,Z^\prime}_{A_t})_{t\geq0}$
and $({\cal F}^{G,Z,Z^\prime}_t)_{t\geq0}$,
respectively.
\end{proposition}

\begin{pf}
First note that $\Gamma$ is a continuous, strictly increasing,
perfect additive functional of $(G,Z)$ which satisfies
(\ref{inftyprop}) and recall that $(G,Z)$ is a strong Markov process
with respect to ${\cal F}^{G,Z,Z^\prime}_{t}, t\ge0$,
by Corollary \ref{basic}. So
$(G\circ A,Z\circ A)=(\tilde{X},\tilde{Y})$
must possess the strong Markov property
with respect to ${\cal F}^{G,Z,Z^\prime}_{A_t}, t\ge0$,
by~\cite{sharpe}, Theorem~65.9.
But $A=\Gamma^{-1}=\int_0^\cdot\tilde{Y}^2_s\,ds$ by time-changing the
integral defining $\Gamma$.
So $\tilde{Y}$ is a continuous-time Markov chain
with $Q$-matrix $(q[y_i,y_j]), y_i,y_j\in\bfS$.
Combining these statements, $(\tilde{X},\tilde{Y})$ has the same law as
$({X},{Y})$ under $P_{x,y}$, since both pairs
satisfy the equation $d{X} = a({X}){Y}\,d{B}$
in the sense explained in the proposition and this equation
admits a weakly unique solution.
As a consequence it follows from (\ref{one prob space}) that
%
%e2.6 #&#
\begin{equation}
\label{with snell} v(x,y)=\sup_{ 0\le\tau\le T } \tilde{E}
\bigl[e^{-q\tau}g(\tilde {X}_{\tau})\bigr],
\end{equation}
where the finite stopping times $\tau$ are with respect to the filtration
${\cal F}^{G,Z,Z^\prime}_{A_t}$, $t\ge0$. Here one should mention that
the stopping times used in (\ref{one prob space}) are with respect to
the filtration generated by $(\tilde{X},\tilde{Y})$ which might be
smaller than ${\cal F}^{G,Z,Z^\prime}_{A_t}, t\ge0$.
However, it is well known
that the corresponding suprema are the same
if the underlying process, in this case $(\tilde{X},\tilde{Y})$,
is also strong Markov with respect to the bigger filtration.
For completeness we sketch the proof of (\ref{with snell})
in the \hyperref[app]{Appendix} on page \pageref{snelli}.

It remains to show that
%
%e2.7 #&#
\begin{equation}
\label{just the same} \sup_{\tau\in\mathcal{T}_T } \tilde{E} \bigl[e^{-q \tau} g(
\tilde{X}_\tau)\bigr]=\sup_{\rho
\in\mathcal{M}_T} \tilde{E}
\bigl[e^{-q\Gamma_\rho} g(G_{\rho})\bigr].
\end{equation}
Fix $\tau\in\mathcal{T}_T$, and observe that
\[
\tilde{E} \bigl[e^{-q \tau}g(\tilde{X}_\tau)\bigr]= \tilde{E}
\bigl[e^{-q\Gamma_{ A_\tau}}g(G_{A_\tau})\bigr]
\]
by property (P1) and the construction of $\tilde{X}$.
Also $A_\tau$ is an
${\cal F}^{G,Z,Z^\prime}_t$ stopping time
by Lemma \ref{lem_equiv_filt}.
The right-hand side above does not change if a finite
version of $A_\tau$ is chosen which still is an
${\cal F}^{G,Z,Z^\prime}_t$ stopping time,
since the filtration satisfies the usual conditions.
Thus $A_\tau\in\mathcal{M}_T$,
and it follows that
\[
\tilde{E} \bigl[e^{-q \tau}g(\tilde{X}_\tau)\bigr]\le \sup
_{\rho
\in\mathcal{M}_T} \tilde{E} \bigl[e^{-q\Gamma_\rho} g(G_{\rho})
\bigr].
\]
Similarly, for fixed $\rho\in\mathcal{M}_T$, the equality
$\tilde{E} [e^{-q\Gamma_\rho}g(G_{\rho})]
=\tilde{E} [e^{-q\Gamma_\rho}g(\tilde{X}_{\Gamma_\rho})]$
leads to
\[
\tilde{E} \bigl[e^{-q\Gamma_\rho}g(G_\rho)\bigr]\le \sup
_{\tau\in\mathcal{T}_T } \tilde{E} \bigl[e^{-q
\tau}g(\tilde{X}_\tau)
\bigr]
\]
finally proving (\ref{just the same}).
\end{pf}

Of course, the conclusion of Proposition \ref{equiv stop probl}
remains valid for
$v(x,y^\prime), \tilde{X}^\prime$, $\tilde{Y}^\prime$,
$\mathcal{T}_T^\prime, \mathcal{M}_T^\prime$,
$A^\prime$
and $\Gamma^\prime$ if these objects are constructed by using
$Z^\prime$ instead of $Z$. Notice that the
solution $G$ is the same.

We are now in the position to formulate and prove the main result
of this section about the validity of (\ref{want to compare}).
The following notion of a skip-free Markov chain is needed:
a continuous-time Markov chain with $Q$-matrix $(q[y_i,y_j])$
taking the states $y_1<\cdots<y_m$ is called \textit{skip-free}\label{skippy}
if the matrix $Q$ is tridiagonal.
%
%th2.5 #&#
\begin{theorem}\label{main one}
Let $(X,Y)$ be a strong Markov process
given on a family of probability spaces
$(\Omega, \mathcal{F}, P_{x,y},(x,y)\in\bfR\times\bfS)$
and let $g\dvtx \bfR\to\bfR$ be a measurable gain function
such that $\{g\ge0\}\neq\varnothing$.
Assume (\ref{eq:integrabilityCond}), that $(X,Y)$ satisfies
conditions \textup{(C1)}, \textup{(C2)} on page \pageref{a and b} and condition
\textup{(C3)} on page \pageref{new C3}
and that all pairs of processes satisfying conditions \textup{(C1)}, \textup{(C2)}
have the same law.
Further suppose that $Y$ is skip-free.
Define ${\cal K}^{g+}_T$ to be the collection
of all finite stopping times $\tau\le T$
with respect to the filtration generated by $(X,Y)$
such that $g(X_\tau)\ge0$.
Fix $(x,y)\in\bfR\times\bfS$ and assume that
$v(x,y) = \sup_{\tau\in{\cal K}^{g+}_T}E_{x,y} [e^{-q\tau}g(X_\tau)]$.
Then
\[
v(x,y)\le v\bigl(x,y^\prime\bigr) \qquad\mbox{for all $y^\prime\in
\bfS$ such that $y\le y^\prime$,}
\]
so that $v(x,y)$ is a lower bound for $v(x,\cdot)$ on $[y,\infty)\cap
\bfS$.
\end{theorem}
%
%re2.6 #&#
\begin{remark}\label{techni}
(i) The condition
$v(x,y) = \sup_{\tau\in{\cal K}^{g+}_T}E_{x,y} [e^{-q\tau}g(X_\tau)]$
is a technical condition
which states that the optimum $v(x,y)$ as defined by (\ref{osp})
can be achieved by stopping at nonnegative values of $g$ only.
It is of course trivially satisfied for all $(x,y)\in\bfR\times\bfS$
if the gain function is nonnegative
and in this case the theorem means that $v(x,\cdot)$ is increasing.

\begin{longlist}[(iii)]
\item[(ii)]In the case of an infinite time horizon $T=\infty$,
it easily follows from the section theorem \cite{Rev-Yor}, Theorem~IV.5.5, that
\[
P_{x,y} \bigl(\inf\bigl\{t\ge0\dvtx
g(X_t)\ge0\bigr\}<\infty \bigr) = 1\qquad
\mbox{for all $(x,y)\in\bfR\times\bfS$}
\]
is sufficient for
$v(x,y) = \sup_{\tau\in{\cal K}^{g+}_T}E_{x,y} [e^{-q\tau}g(X_\tau)]$
to be true for all $(x,y)\in\bfR\times\bfS$ since $(X,Y)$ is strong Markov.
Indeed, if a process always hits the set $\{g\ge0\}$ with probability
one, then it is quite natural that maximal gain is obtained while avoiding
stopping at negative values of $g$.
One can easily construct processes satisfying this sufficient condition
where the gain function $g$ takes both positive
\emph{and} negative values.

\item[(iii)]In the case where $T<\infty$,
the only reasonable sufficient condition the authors can find is the
trivial condition $g(x)\ge0$ for all $x\in\bfR$.
This is because, in general a process is not almost surely guaranteed
to hit a subset of the state space in finite time.

\item[(iv)]
The monotonicity result of this theorem supports the intuition that
the larger the diffusion coefficient (volatility) of a diffusion
without drift,
the faster this diffusion moves and hence the sooner it reaches the points
where the gain function $g$ is large. As the killing term of the cost functional
defining $v(x,y)$ punishes the elapsed time, $v(x,y')$ should indeed be larger
than $v(x,y)$, for $y'>y$, if the volatility process starting from $y'$
stays above
the volatility process starting from $y$, and this is ensured by the
skip-free property
of the Markov chain.
\end{longlist}
\end{remark}
\begin{pf*}{Proof of Theorem \ref{main one}}
Fix $x\in\bfR$ and $y,y^\prime\in\bfS$ such that $y\le y^\prime$,
and let $G,W,{Z},{Z}^\prime$ be given on a complete probability space
$(\tilde{\Omega}, \tilde{\mathcal{F}},\tilde{P})$ as described in
Corollary \ref{basic}. While in Corollary \ref{basic} the coupling of
the two chains $Z$ and $Z^\prime$ was not specified any further we now
choose a particular coupling
associated with a $Q$-matrix $\mathfrak{Q}$
which allows us to compare $Z$ and $Z^\prime$ directly.
Denoting the $Q$-matrix corresponding to the independence
coupling by $\mathfrak{Q}^\perp$, we set
\[
\mathfrak{Q}
\left[\begin{array}{c@{\hspace*{4pt}}|@{\hspace*{4pt}}c}
y_i&y_j\\y_k&y_l
\end{array}\right]
%y_k&y_l }\right]
 = \cases{
\mathfrak{Q}^\perp
%y_k&y_l }
\left[\begin{array}{c@{\hspace*{4pt}}|@{\hspace*{4pt}}c}
y_i&y_j\\y_k&y_l
\end{array}\right]
,&\quad $i\neq k,$\vspace*{2pt}\cr
y_i^{-2} q[y_i,y_j],&\quad
$i=k, j=l,$\vspace*{2pt}\cr
0,&\quad $i=k, j\neq k$ }
\]
for $y_i,y_j,y_k,y_l\in\bfS$; that is,
$Z$ and $Z^\prime$ move independently until they
hit each other for the first time and then they move together.
It follows from the skip-free-assumption that
$Z$ cannot overtake $Z^\prime$
before they hit each other for the first time.
Hence
\[
Z_0=y\le y^\prime=Z_0^\prime\qquad
\mbox{implies } Z_t\le Z_t^\prime, t\ge0,
\mbox{a.s.},
\]
which results in the inequality
%
%e2.8 #&#
\begin{equation}
\label{fundamental compari} \Gamma_t =\int_0^t
Z_s^{-2}\,ds \ge \int_0^t
\bigl(Z_s^\prime\bigr)^{-2}\,ds =
\Gamma^\prime_t,\qquad t\ge0, \mbox{a.s.}
\end{equation}
Note that then the inverse increasing processes
$A=\Gamma^{-1}$ and $A^\prime=(\Gamma^\prime)^{-1}$
must satisfy the relation $A_t\le A^\prime_t, t\ge0$, a.s.

Now recall the definition of $\mathcal{M}_T$
in Proposition \ref{equiv stop probl},
and note that the above comparison allows us to conclude that
%
%e2.9 #&#
\begin{equation}
\label{compari} \tilde{E} \bigl[e^{-q\Gamma_\rho}g(G_\rho)\bigr] \le
\tilde{E} \bigl[e^{-q\Gamma^\prime_\rho}g(G_\rho)\bigr] \qquad\mbox{for every $\rho
\in\mathcal{M}_T^+$},
\end{equation}
where
$\mathcal{M}_T^+=
\{\rho\in\mathcal{M}_T\dvtx g(G_\rho)\ge0 \mbox{ a.s.}\}$.
Thus
\[
\tilde{E} \bigl[e^{-q\Gamma_\rho}g(G_\rho)\bigr] \le \sup
_{\rho^\prime\in\mathcal{M}^\prime_T} \tilde{E} \bigl[e^{-q\Gamma^\prime_{ \rho^\prime}} g(G_{\rho^\prime})
\bigr] \qquad\mbox{for every $\rho\in\mathcal{M}_T^+$}
\]
since $A_T\le A_T^\prime$ a.s. implies that
every stopping time in $\mathcal{M}_T^+$ has a version which is
in $\mathcal{M}^\prime_T$. Putting these results together, we obtain
\[
\sup_{\rho\in\mathcal{M}_T^+} \tilde{E} \bigl[e^{-q\Gamma_{ \rho}}
g(G_{\rho})\bigr] \le \sup_{\rho^\prime\in\mathcal{M}^\prime_T} \tilde{E}
\bigl[e^{-q\Gamma^\prime_{ \rho^\prime}} g(G_{\rho^\prime})\bigr].
\]
But, if $\mathcal{T}_T^+$ denotes
$\{\tau\in\mathcal{T}_T\dvtx g(\tilde{X}_\tau)\ge0 \mbox{ a.s.}\}$,
then the equality
\[
\sup_{\tau\in\mathcal{T}_T^+ } \tilde{E} \bigl[e^{-q \tau} g(
\tilde{X}_\tau)\bigr]=\sup_{\rho
\in\mathcal{M}_T^+} \tilde{E}
\bigl[e^{-q\Gamma_\rho} g(G_{\rho})\bigr]
\]
can be shown in the same way that (\ref{just the same}) was shown in
the proof
of Proposition \ref{equiv stop probl} (note that in this proof
we may choose versions of certain stopping times and this is the reason
the qualification ``a.s.'' appears in the definitions of $\mathcal
{M}_T^+$ and
$\mathcal{T}_T^+$).

Furthermore,
\[
\mbox{if } v(x,y) = \sup_{\tau\in{\cal K}^{g+}_T}E_{x,y} \bigl[
e^{-q\tau}g(X_\tau)\bigr] \mbox{ then } v(x,y) \le \sup
_{\tau\in\mathcal{T}_T^+ } \tilde{E} \bigl[e^{-q \tau} g(\tilde{X}_\tau)
\bigr]
\]
since the law of $(\tilde{X},\tilde{Y})$ is equal to the law of
$({X},{Y})$ under $P_{x,y}$ and the filtration
${\cal F}^{G,Z,Z^\prime}_{A_t}, t\ge0$, is at least as big as the
filtration generated by $(\tilde{X},\tilde{Y})$.
So, under the condition
$v(x,y) =
\sup_{\tau\in{\cal K}^{g+}_T}E_{x,y} [ e^{-q\tau}g(X_\tau)]$,
we can finally deduce that
\begin{eqnarray*}
v(x,y)&\le& \sup_{\tau\in\mathcal{T}_T^+ } \tilde{E} \bigl[e^{-q \tau} g(
\tilde{X}_\tau)\bigr]
\\
&=& \sup_{\rho\in\mathcal{M}_T^+} \tilde{E} \bigl[e^{-q\Gamma_\rho}g(G_\rho)
\bigr] \le \sup_{\rho^\prime\in\mathcal{M}_T^\prime} \tilde{E} \bigl[e^{-q\Gamma^\prime_{ \rho^\prime}}g(G_{\rho^\prime})
\bigr] = v\bigl(x,y^\prime\bigr),
\end{eqnarray*}
where the last equality is due to Proposition \ref{equiv stop probl}
applied to $(\tilde{X}^\prime,\tilde{Y}^\prime)$.
\end{pf*}
%
%co2.7 #&#
\begin{corollary}\label{version for negative}
\textup{(i)} If $\{g\ge0\}=\varnothing$,
but all other assumptions of Theorem~\ref{main one}
are satisfied, then in the infinite time horizon case where $T=\infty$,
\[
v(x,y) \ge v\bigl(x,y^\prime\bigr) \qquad\mbox{for all $x\in\bfR$ and
$y,y^\prime\in\bfS$ such that $y\le y^\prime$,}
\]
so that $v(x,\cdot)$ is decreasing.

\begin{longlist}[(ii)]
\item[(ii)] Let the assumptions of
Theorem \ref{main one} be based on $\bfS\subseteq(-\infty,0)$,
fix $(x,y)\in\bfR\times\bfS$ and assume that
$v(x,y) = \sup_{\tau\in{\cal K}^{g+}_T}E_{x,y} [e^{-q\tau}g(X_\tau)]$.
Then
\[
v(x,y)\ge v\bigl(x,y^\prime\bigr) \qquad\mbox{for all $y^\prime\in
\bfS$ such that $y\le y^\prime$,}
\]
so that $v(x,y)$ is an upper bound for $v(x,\cdot)$ on
$[y,\infty)\cap\bfS$.
\end{longlist}
\end{corollary}
\begin{pf}
If $\{g\ge0\}=\varnothing$ then, instead of (\ref{compari}), we obtain
\[
\tilde{E} \bigl[e^{-q\Gamma_\rho}g(G_\rho)\bigr] \ge \tilde{E}
\bigl[e^{-q\Gamma^\prime_\rho}g(G_\rho)\bigr]\qquad \mbox{for every $\rho\in
\mathcal{M}_T$}
\]
and $\mathcal{M}_T$ is (up to versions) equal to $\mathcal
{M}_T^\prime$
since $T=\infty$.
Hence (i) can be deduced directly from Proposition \ref{equiv stop probl}.
Note that the above inequality cannot be used in the case where
$T<\infty$ since there can be stopping times in $\mathcal{M}_T^\prime$
which are not in $\mathcal{M}_T$.

If $\bfS\subseteq(-\infty,0)$ then $Z_t\le Z_t^\prime, t\ge0$,
a.s., does
not imply (\ref{fundamental compari}) but instead
\[
\Gamma_t =\int_0^t
Z_s^{-2}\,ds \le \int_0^t
\bigl(Z_s^\prime\bigr)^{-2}\,ds =
\Gamma^\prime_t,\qquad t\ge0, \mbox{a.s.},
\]
hence, interchanging the roles of $y$ and $y^\prime$,
(ii) can be proved like Theorem~\ref{main one}.
Note that $Z_t\le Z_t^\prime, t\ge0$, a.s., would not lead to any
comparison between $\Gamma$ and $\Gamma^\prime$ if $y<0<y^\prime$.
Hence our method cannot be applied to show the monotonicity of
$v(x,y)$ in $y\in\bfS$ if $\bfS$ contains a neighbourhood of zero.
\end{pf}

%s3 #&#
\section{The diffusion case}\label{diffu}

Fix $\delta\in[-1,1]$, and
suppose that $({X},{Y})$ is a strong Markov process
given on a family of probability spaces
$({\Omega}, {\mathcal{F}},P_{x,y}$, $(x,y)\in\bfR\times\bfS)$
which satisfies the following conditions:
\begin{longlist}[(C1$'$)]\label{a,b and c}
\item[(C1$'$)]
the process $(X,Y)$ is adapted with respect to a filtration
${\mathcal{F}}_t, t\ge0$, of sub-$\sigma$-algebras
of ${\mathcal{F}}$ and,
for every $(x,y)\in\bfR\times\bfS$, there is a pair $({B},{B}^Y)$ of
${\mathcal{F}}_t$ Brownian motions
on $({\Omega}, {\mathcal{F}}, P_{x,y})$ with covariation
$\langle{B}, {B}^Y \rangle_t=\delta t, t\ge0$, such that
\[
{X}_t = {X}_0+\int_0^t
a({X}_s){Y}_s \,d{B}_s \quad\mbox{and}\quad
{Y}_t = {Y}_0+\int_0^t
\eta({Y}_s)\,d{B}_s^Y + \int
_0^t\theta({Y}_s)\,ds
\]
for all $t\ge0, P_{x,y}$-a.s.;
\item[(C2$'$)]
the process $Y$ takes values in $\bfS\subseteq(0,\infty)$ and
\[
P_{x,y} \biggl( \lim_{t\uparrow\infty} \int
_0^t {Y}^2_s \,ds = \infty
\biggr)=1 \qquad\mbox{for all } (x,y)\in\bfR\times\bfS.
\]
\end{longlist}
%
%re3.1 #&#
\begin{remark}\label{weak solution}
Under the assumptions above,
for every $(x,y)\in\bfR\times\bfS$,
there exists a weak solution to the system of
stochastic differential equations
%
%e3.1 #&#
\begin{equation}
\label{equ star} \cases{ %
dG = a(G) \,dW,
\vspace*{2pt}\cr
d\xi= {\eta(\xi)} {\xi^{-1}} \,dW^\xi+ {
\theta(\xi)} {\xi^{-2}} \,dt,
\vspace*{2pt}\cr
\xi_t\in\bfS, \qquad t\ge0,}
\end{equation}
driven by a pair of Brownian motions with covariation
$\langle{W}, {W}^\xi\rangle_t=\delta t, t\ge0$.
Such a solution can be given by
$(X\circ\langle M\rangle^{-1},Y\circ\langle M\rangle^{-1})$
where $M$ denotes the continuous local martingale
$M_s =\int_0^s Y_u \,dB_u, s\ge0$,
as in Remark \ref{G strong markov}(ii).
Here $W=M\circ\langle M\rangle^{-1}$ and
${W}^\xi=\int_0^{\langle M\rangle^{-1}_\cdot}Y_s \,dB^Y_s$
are ${\mathcal{F}}_{\langle M\rangle^{-1}_t}$\vspace*{1pt} Brownian motions
by Dambis--Dubins--Schwarz' theorem (see \cite{Rev-Yor}, Theorem~V.1.6)
with covariation
\begin{eqnarray*}
\bigl\langle W, W^\xi\bigr\rangle_t&=& \biggl\langle\int
_0^\cdot Y_s \,dB_s, \int
_0^{\cdot}Y_s \,dB^Y_s
\biggr\rangle_{\langle M\rangle^{-1}_t}
\\
&=& \delta\int_0^{\langle M\rangle^{-1}_t}Y_s^{2}\,ds
= \delta\langle M\rangle_{\langle M\rangle^{-1}_t} = \delta t,\qquad t\ge0,
P_{x,y}\mbox{-a.s.},
\end{eqnarray*}
where the last equality is ensured by condition (C2$'$).
\end{remark}

We want to show (\ref{want to compare}) using a method similar
to the method applied in Section~\ref{regswi}.
The main difference to the case discussed in Section~\ref{regswi} is
that the pair $({X},{Y})$ is now determined by a system of
stochastic differential equations.
So, instead of constructing $\tilde{X}$
by time-changing a solution of the single
equation $dG=a(G)\,dW$ as in Section~\ref{regswi},
we now construct $\tilde{X}$ by time-changing a solution
of a system of stochastic differential equations.
Furthermore, in Section~\ref{regswi} we constructed
the coupling of $Z$ and $Z^\prime$ in the proof of Theorem \ref{main one}
from a given generator.
In this section we will couple $\xi$ and $\xi^\prime$---both
satisfying the second equation in (\ref{equ star}) but starting from
$y\le y^\prime$, respectively---we will do so directly from
the stochastic differential equation. As a consequence, the next condition
appears to be slightly stronger than the corresponding condition (C3) of the last section.
However, in Theorem~\ref{main one} we needed (C3), a skip-free
Markov chain
and weak uniqueness of~(\ref{eq:dynX}) while below, in the
corresponding Theorem \ref{main two},
we will only need:
\begin{longlist}[(C3$'$)]
\item[(C3$'$)]\label{only c}
Let $a,\eta,\theta$ be measurable functions such that the
system of stochastic differential equations
(\ref{equ star})
has, for all initial conditions $(G_0,\xi_0)\in\bfR\times\bfS$,
a unique nonexploding strong solution
taking values in $\bfR\times\bfS$.
\end{longlist}

Now choose a complete probability space
$(\tilde{\Omega}, \tilde{\mathcal{F}}, \tilde{P})$
big enough to carry a pair of Brownian motions $(W,W^\xi)$
with covariation
$\langle{W}, {W}^\xi\rangle_t=\delta t, t\ge0$, and denote by
$\tilde{\mathcal{F}}_t, t\geq0$,
the usual augmentation of the filtration generated by $(W,W^\xi)$.
Let $(G,\xi)$ be the unique solution of the system (\ref{equ star})
starting from $G_0=x\in\bfR$ and $\xi_0=y\in\bfS$
given on $(\tilde{\Omega}, \tilde{\mathcal{F}}, \tilde{P})$
by $(W,W^\xi)$.

Define $\Gamma=(\Gamma_t)_{t\ge0}$ by
\[
\Gamma_t = \int_0^t{
\xi_u^{-2}}\,du,\qquad t\ge0,
\]
and remark that $\Gamma$ satisfies (\ref{inftyprop}). Indeed,
by Remark \ref{weak solution}, $Y\circ\langle M\rangle^{-1}$
solves the second equation of (\ref{equ star}), and hence condition
(C3$'$) implies
that $\xi$ has the same law as $Y\circ\langle M\rangle^{-1}$
under $P_{x,y}$.
Property (\ref{inftyprop}) therefore follows from (C2$'$) since
$\int_0^\cdot{Y}_{\langle M\rangle^{-1}_s}^{-2} \,ds$
is the inverse of $\int_0^\cdot Y_s^2\,ds$.

Of course, we may deduce from
(\ref{inftyprop}) together with the fact that $\xi$ never vanishes,
that
$\Gamma$ is a continuous and strictly increasing process.
Thus, $A=\Gamma^{-1}$
is also a continuous and strictly increasing process
satisfying (\ref{inftyprop for A}). As a consequence, the two technical
properties (P1) and (P2) on page \pageref{prop_tech} must
again be valid.

As $\xi$ is $\tilde{\mathcal{F}}_t$-adapted, we see that
%
%e3.2 #&#
%e3.3 #&#
\begin{eqnarray}
G_t&=& x+\int_0^t
a(G_s)\xi_s \,d\tilde{M}_s,\qquad t\ge0,
\mbox{a.s.},\label{make G}
\\
\xi_t&=& y+\int_0^t\eta(
\xi_s)\,d\tilde{M}^\xi_s+\int
_0^t\theta(\xi _s)\,d
\Gamma_s > 0,\qquad t\ge0, \mbox{a.s.},\label{make xi}
\end{eqnarray}
where (\ref{inftyprop}) implies that the continuous local martingales
$\tilde{M}$ and $\tilde{M}^\xi$ given by the
stochastic integrals
\[
\tilde{M}_s =\int_0^s
\xi_u^{-1}\,dW_u \quad\mbox{and}\quad
\tilde{M}^\xi_s =\int_0^s
\xi_u^{-1}\,dW^\xi_u
\]
exist for each $s\ge0$.
Now it immediately follows from (\ref{make G}), (\ref{make xi}) that
the $\tilde{\mathcal{F}}_{A_t}$-adapted processes
$\tilde{X}=G\circ A$ and $\tilde{Y}=\xi\circ A$
on $(\tilde{\Omega}, \tilde{\mathcal{F}}, \tilde{P})$
constitute a nonexploding weak solution
of the system (\ref{eq:dynX}), (\ref{vol_W}) with
$\tilde{Y}_t\in\bfS, t\ge0$,
since $\tilde{B}=\tilde{M}\circ A$ and $\tilde{B}^Y=\tilde{M}^\xi
\circ A$
are $\tilde{\mathcal{F}}_{A_t}$ Brownian motions
by Dambis--Dubins--Schwarz' theorem \cite{Rev-Yor}, Theorem~V.1.6 and
\[
\bigl\langle\tilde{B}, \tilde{B}^Y \bigr\rangle_t= \bigl
\langle\tilde{M}, \tilde{M}^\xi\bigr\rangle_{A_t}= \delta
\int_0^{A_t}\xi_u^{-2}\,du=
\delta\Gamma_{A_t}=\delta t,\quad t\ge0, \mbox{a.s.},
\]
by property (P1).
%
%re3.2 #&#
\begin{remark}\label{reverting}
(i) Combining Remark \ref{weak solution} and condition (C3$'$),
it follows from
the construction above that $(\tilde{X},\tilde{Y})$
must have the same distribution as $({X},{Y})$ under $P_{x,y}$.
\begin{longlist}[(ii)]
\item[(ii)] The filtration $\tilde{\mathcal{F}}_{A_t}, t\ge0$,
might be
bigger than the filtration generated by $(\tilde{X},\tilde{Y})$.
However, it is straightforward to show the strong Markov property of
$(\tilde{X},\tilde{Y})$ with respect to
$\tilde{\mathcal{F}}_{A_t}, t\ge0$, since $(\tilde{X},\tilde{Y})$
was obtained by time-changing a unique strong solution of a
system of stochastic differential equation driven by Brownian motions.
\end{longlist}
\end{remark}
This remark makes clear that the following proposition
can be proved by applying the ideas used in the proof of
Proposition \ref{equiv stop probl}
in Section~\ref{regswi} (so we omit its proof).
%
%pr3.3 #&#
\begin{proposition}\label{equiv stop probl for diffu}
Let $G,\tilde{X},\tilde{Y}$
be the processes on the filtered probability space
$(\tilde{\Omega}, \tilde{\mathcal{F}}, \tilde{\mathcal{F}}_t$,
$ t\ge0, \tilde{P})$
introduced above and starting from
$G_0=\tilde{X}_0=x\in\bfR$ and $\tilde{Y}_0=y\in\bfS$.
Then, for any $T\in[0,\infty]$, it follows that
\[
v(x,y)=\sup_{\tau\in\mathcal{T}_T } \tilde{E} \bigl[e^{-q \tau} g(
\tilde{X}_\tau)\bigr]=\sup_{\rho
\in\mathcal{M}_T} \tilde{E}
\bigl[e^{-q\Gamma_\rho} g(G_{\rho})\bigr],
\]
where
\[
\mathcal{T}_T=\{\tau\in\mathcal{T} \dvtx 0\leq\tau\leq T\}
\quad\mbox{and} \quad\mathcal{M}_T=\{\rho\in\mathcal{M} \dvtx 0\leq\rho\leq
A_T\}.
\]
Here,
$\mathcal{T}$ and $\mathcal{M}$
denote the families of finite stopping times with respect to the
filtrations
$(\tilde{\mathcal{F}}_{A_t})_{t\geq0}$
and $(\tilde{\mathcal{F}}_t)_{t\geq0}$,
respectively.
\end{proposition}
%
%re3.4 #&#
\begin{remark}
The above representation of %the value function
$v(x,y), (x,y)\in\bfR\times\bfS$, could be extended to cases where
$\bfS$ is bigger than $(0,\infty)$. However, in such cases,
the equation for $\xi$ in (\ref{equ star})
must admit solutions starting from $\xi_0=0$ which is an additional constraint,
since $\xi$ is in the denominator on the right-hand side of this
equation. Furthermore, in addition to the assumption that
$P_{x,y} ( \lim_{t\uparrow\infty} \int_0^t {Y}^2_s \,ds =\infty)=1$
one would need to assume that $\int_0^\cdot{Y}^2_s \,ds$ is
strictly increasing $P_{x,y}$-a.s. as, in principle,
the process $Y$ could now spend time at zero.
\end{remark}

Recall that, in contrast to the case of regime-switching, the
process $\tilde{Y}$ above was constructed by time-change from a
solution of a
stochastic differential equation and this results in some small
variations from the proof of Theorem \ref{main one}.
Note that the conclusion of Proposition \ref{equiv stop probl for
diffu} remains valid for
$v(x,y^\prime), \tilde{X}^\prime, \mathcal{T}_T^\prime, \mathcal
{M}_T^\prime, A^\prime$
and $\Gamma^\prime$ if these objects are constructed using a
different starting point $y^\prime\in\bfS$.
%
%th3.5 #&#
\begin{theorem}\label{main two}
Let $(X,Y)$ be a strong Markov process
given on a family of probability spaces
$(\Omega, \mathcal{F}, P_{x,y},(x,y)\in\bfR\times\bfS)$,
and let $g\dvtx \bfR\to\bfR$ be a measurable gain function
such that $\{g\ge0\}\neq\varnothing$.
Assume (\ref{eq:integrabilityCond}), that $(X,Y)$ satisfies
conditions \textup{(C1$'$)} and \textup{(C2$'$)} on page \pageref{a,b and c}
and that
condition \textup{(C3$'$)} on page \pageref{only c} holds true for system
(\ref{equ star}).
Define ${\cal K}^{g+}_T$ to be the collection
of all finite stopping times $\tau\le T$
with respect to the filtration generated by $(X,Y)$
such that $g(X_\tau)\ge0$.
Fix $(x,y)\in\bfR\times\bfS$ and assume that
$v(x,y) = \sup_{\tau\in{\cal
K}^{g+}_T}E_{x,y} [e^{-q\tau}g(X_\tau)]$.
Then
\[
v(x,y)\le v\bigl(x,y^\prime\bigr)\qquad \mbox{for all $y^\prime\in
\bfS$ such that $y\le y^\prime$,}
\]
so that $v(x,y)$ is a lower bound for $v(x,\cdot)$ on
$[y,\infty)\cap\bfS$.
\end{theorem}

\begin{pf}
Fix $x\in\bfR$ and $y,y^\prime\in\bfS$
with $y\le y^\prime$
and choose a complete probability space
$(\tilde{\Omega}, \tilde{\mathcal{F}}, \tilde{P})$
large enough to carry a pair of Brownian motions $(W,W^\xi)$
with covariation
$\langle{W}, {W}^\xi\rangle_t=\delta t, t\ge0$.
Let $(G,\xi)$ and $(G,\xi^\prime)$ be the solutions of~(\ref{equ star})
starting from $(x,y)$ and from $(x,y^\prime)$,
respectively, which are both given by $(W,W^\xi)$
on $(\tilde{\Omega}, \tilde{\mathcal{F}}, \tilde{P})$.
Remark that $G$ is indeed the same for both pairs since~(\ref{equ
star}) is a
system of decoupled equations.

Define $C=\inf\{t\ge0\dvtx \xi_t>\xi^\prime_t\}$ and set
$\bar{\xi}_t=\xi_t I(t<C)+\xi^\prime_t I(t\ge C)$
so that $\bar{\xi}_t\le\xi^\prime_t$ for all $t\ge0$.
Obviously, $(G,\bar{\xi})$ solves system (\ref{equ star}) starting from
$G_0=x$ and $\bar{\xi}_0=y$ hence
$\xi_t=\bar{\xi}_t\le\xi^\prime_t, t\ge0$, a.s., by strong
uniqueness.

Construct $\tilde{X}=G\circ A$ and $\tilde{X}^\prime=G\circ A^\prime$
using the above $(G,\xi)$ and $(G,\xi^\prime)$,
and observe that
$\Gamma_t\geq\Gamma^\prime_t$ follows immediately
from $0<\xi_t\le\xi^\prime_t$ for all $t\ge0$ a.s.
Thus, simply using Proposition \ref{equiv stop probl for diffu}
instead of Proposition \ref{equiv stop probl},
the rest of the proof can be copied from the corresponding part
of the proof of Theorem~\ref{main one}.
\end{pf}
%
%re3.6 #&#
\begin{remark}
For a discussion of the technical condition
$v(x,y) =\break   \sup_{\tau\in{\cal K}^{g+}_T}E_{x,y} [e^{-q\tau}g(X_\tau)]$
we refer the reader to Remark \ref{techni}.
Corollary \ref{version for negative} remains true if it
is reformulated in the terms of the theorem above
instead of Theorem \ref{main one}.
\end{remark}

%s4 #&#
\section{Continuity in the diffusion-case}\label{conti}

Let $\bfS$ be an open subset of $(0,\infty)$,
fix $x\in\bfR$
and suppose that all the assumptions of Theorem \ref{main two} are satisfied.
Furthermore, suppose that
%
%e4.1 #&#
\begin{equation}
\label{Kg for sequence} v(x,y) = \sup_{\tau\in{\cal K}^{g+}_T}E_{x,y}
\bigl[e^{-q\tau}g(X_\tau)\bigr]\qquad \mbox{for all $y\in\bfS$.}
\end{equation}
For a sequence $(y_n)_{n=1}^\infty\subseteq\bfS$
converging to $y_0\in\bfS$ as $n\to\infty$,
denote by $(G,\xi^n)$ the solution of (\ref{equ star}),
starting from $G_0=x$ and $\xi^n_0=y_n, n=0,1,2,\ldots,$
given by a pair
$(W,W^\xi)$ of Brownian motions with covariation
$\langle{W}, {W}^\xi\rangle_t=\delta t, t\ge0$,
on a probability space
$(\tilde{\Omega}, \tilde{\mathcal{F}}, \tilde{P})$.
Using $(G,\xi^n)$ construct $\Gamma^n,A^n, n=0,1,2,\ldots,$
like~$\Gamma,A$ in Section~\ref{diffu}.
%
%le4.1 #&#
\begin{lemma}\label{with feller}
Suppose that the $\xi$-component corresponding to the unique strong solution
to (\ref{equ star}) on page \pageref{equ star}
is a Feller process with state space $\bfS$.
If the sequence $(y_n)_{n=1}^\infty$ is monotone, that is,
either $y_n\downarrow y_0$ or $y_n\uparrow y_0$ when $n\to\infty$,
then
%
%e4.2 #&#
\begin{equation}
\label{time ch conv} \Gamma^n_t\to\Gamma^0_t
\quad\mbox{and}\quad A^n_t\to A^0_t\qquad \mbox{as } n\to
\infty\mbox{ for all } t\ge0 \mbox{ a.s.}
\end{equation}
\end{lemma}

\begin{pf}
Here, we will use without further comment the elementary fact that if
$U$ and $V$ are two random variables with the same law and $U\ge V$
a.s. then, in fact, $U=V$ a.s.

Suppose that $y_n\downarrow y_0$ as $n\to\infty$.
By the coupling argument in the proof of Theorem \ref{main two},
without loss of generality one may chose $\xi^n$ such that
%
%e4.3 #&#
\begin{equation}
\label{consec_couplings} \xi^1_t \geq\xi_t^2
\geq\cdots\geq\xi_t^n \geq\cdots \geq
\xi^0_t > 0,\qquad t\ge0, n=0,1,2,\ldots;
\end{equation}
hence the pathwise limit $\lim_n\xi^n_t, t\ge0$, exists.
It follows from the Feller property that the two
processes $\xi^0$ and $(\lim_n\xi^n_t)_{t\ge0}$
must have the same law by comparing their finite-dimensional
distributions.
%but limit only upper semicountinuous
As (\ref{consec_couplings}) also yields the inequalities
$\lim_n\xi^n_t\ge\xi^0_t>0, t\ge0$, we see that
\[
\Gamma^0_t \ge \int_0^t
\Bigl(\lim_n\xi^n_u
\Bigr)^{-2}\,du = \lim_n\Gamma^n_t,\qquad
t\ge0,
\]
by monotone convergence.
But, if $\xi^0$ and $(\lim_n\xi^n_t)_{t\ge0}$
have the same law,
then the same must hold true for
$\Gamma^0$ and $\int_0^t(\lim_n\xi^n_u)^{-2}\,du, t\ge0$.
Thus $\Gamma^0_t \ge\int_0^t(\lim_n\xi^n_u)^{-2}\,du$ implies
$\Gamma^0_t = \int_0^t(\lim_n\xi^n_u)^{-2}\,du$ a.s. for each $t\ge0$. The desired result,
$\Gamma^0_t = \int_0^t(\lim_n\xi^n_u)^{-2}\,du, t\ge0$, a.s., now
follows since
both processes have continuous paths.

Thus
$\Gamma^n_t\uparrow\Gamma^0_t, t\ge0$, a.s.
Since $A^n$,
$A^0$ are the right-inverses of the continuous increasing processes
$\Gamma^n$ and $\Gamma^0$, respectively, we have
$A^n_t\downarrow A^0_t, t\ge0$, a.s., completing the proof in the
case where the $(y_n)$ are decreasing.

In the case where $y_n\uparrow y_0$ as $n\to\infty$, we see that
\[
0 < \xi^1_t \leq\xi_t^2 \leq
\cdots\leq\xi_t^n \leq\cdots\leq \xi^0_t
\]
and
\[
\Gamma^0_t \le \int_0^t
\Bigl(\lim_n\xi^n_u
\Bigr)^{-2}\,du = \lim_n\Gamma^n_t,\qquad
t\ge0,
\]
by Lebesgue's dominated convergence theorem.
This ensures that $\Gamma^n_t\downarrow\Gamma^0_t$, $t\ge0$, a.s., and
$A^n_t\uparrow A^0_t, t\ge0$, a.s.
\end{pf}

In what follows,
in addition to the assumptions of Theorem \ref{main two},
we impose the assumption of Lemma \ref{with feller} and
the following condition (C4$'$)
is used to summarize these conditions, that is:
\begin{itemize}[(C4$'$)]\label{summarising}
\item[(C4$'$)]
\begin{itemize}
\item
the gain function $g$ satisfies $\{g\ge0\}\neq\varnothing$;
\item
the process $(X,Y)$ satisfies conditions
(\ref{eq:integrabilityCond}), (C1$'$), (C2$'$) and the value
function $v$ satisfies (\ref{Kg for sequence}) for the chosen value of $x$;
\item
condition (C3$'$) holds true for the system (\ref{equ star}) and the
second equation in~(\ref{equ star}) has a Feller solution.
\end{itemize}
\end{itemize}
Note that, in many cases, the conditions one imposes on the
coefficients $\eta$ and $\theta$ to ensure
condition (C3$'$) also imply that the whole solution
of (\ref{equ star}) is a Feller process.

We now discuss the continuity of the value function
$v(x,\cdot)$ which we subdivide into left-continuity and right-continuity.

%pr4.2 #&#
\begin{proposition}\label{left conti for inf}
Assume condition \textup{(C4$'$)}. Then, when $T=\infty$,
$v(x,\cdot)$ is left-continuous.
\end{proposition}
\begin{pf}
First observe that Theorem \ref{main two} implies that
\[
\limsup_{n\rightarrow\infty} v(x,y_n) \leq v(x,y_0),
\]
whenever $y_n\uparrow y_0$ in $\bfS$,
so it remains to show that
\[
v(x,y_0) \le\liminf_{n\rightarrow\infty} v(x,y_n).
\]
Recall the definition of $\mathcal{M}$ from Proposition \ref{equiv stop
probl for diffu},
and choose $\rho\in\mathcal{M}$. Then
%
%e4.4 #&#
\begin{equation}
\label{ineq} e^{-q\Gamma^n_\rho}g(G_{\rho}) \ge -e^{-q\Gamma^0_\rho}\bigl
|g(G_{\rho})\bigr|
= -e^{-q\Gamma^0_\rho}\bigl|g\bigl(\tilde{X}^0_{\Gamma^0_\rho}\bigr)\bigr|
\end{equation}
for all $n=1,2,\ldots,$ since $\Gamma^n_\rho\ge\Gamma^0_\rho$.
But the right-hand side of (\ref{ineq}) is integrable by
(\ref{eq:integrabilityCond}). Thus the inequality
%
%e4.5 #&#
\begin{equation}
\label{fatou} \tilde{E} e^{-q\Gamma^0_\rho} g(G_\rho) \leq \liminf
_{n\rightarrow\infty} \tilde{E} e^{-q\Gamma^n_\rho} g(G_\rho)
\end{equation}
follows from Fatou's lemma and Lemma \ref{with feller}.\vspace*{1pt}

Now
$\tilde{E} e^{-q\Gamma^n_\rho}g(G_\rho)\le
\sup_{\rho^\prime\in\mathcal{M}}
\tilde{E} e^{-q\Gamma^n_{ \rho^\prime}}g(G_{\rho^\prime})$,
and so Proposition \ref{equiv stop probl for diffu} gives
%
%e4.6 #&#
\begin{equation}
\label{ineq2} \tilde{E} e^{-q\Gamma^0_\rho} g(G_\rho) \leq \liminf
_{n\rightarrow\infty}v(x,y_n)
\end{equation}
since ${\cal M}_T$ can be replaced by ${\cal M}$ in the case where
$T=\infty$.
So, taking the supremum over $\rho\in\mathcal{M}$ in the
left-hand side of (\ref{ineq2}) completes the proof.
\end{pf}
%
%re4.3 #&#
\begin{remark}\label{lower semi cont}
(i) The fact that $v(x,y_0) \le\liminf_{n\rightarrow\infty} v(x,y_n)$
when $y_n\downarrow y_0$ in $\bfS$ is an immediate
consequence of Theorem \ref{main two}.
As $v(x,y_0)\le\break \liminf_{n\rightarrow\infty} v(x, y_n), y_n\uparrow y_0$,
was shown in the proof above,\vadjust{\goodbreak}
$v(x,\cdot)$ is,
under condition (C4$'$),
lower semicontinuous on $\bfS$
when $T=\infty$ without any continuity-assumption
on the gain function $g$.
\begin{longlist}[(iii)]
\item[(ii)]
From (i) above it follows that,
to establish right-continuity in the case where $T=\infty$,
it remains to show that
$\limsup_{n\rightarrow\infty} v(x,y_n) \leq v(x,y_0)$
when $y_n\downarrow y_0$ in $\bfS$.
We are only able to prove this using the extra
integrability condition of Proposition \ref{new int cond} below.
Note that the combination of Propositions \ref{left conti for inf}
and~\ref{new int cond} gives continuity of
$v(x,\cdot)$ for fixed $x$ in the case where $T=\infty$
without the requirement that the gain function $g$
is continuous.

\item[(iii)]If $T<\infty$, then the proof of Proposition \ref{left
conti for inf} fails.
Indeed, in this case, $\rho$~cannot be chosen from $\mathcal{M}$
as it belongs to a different class $\mathcal{M}^n_T$ for each
$n=0,1,2,\ldots$ We are able to show
left- and right-continuity in the case where $T<\infty$
under the additional assumption that the gain function $g$ is
continuous; see Proposition \ref{conti finite time}.
\end{longlist}
\end{remark}
%
%pr4.4 #&#
\begin{proposition}\label{new int cond}
Assume, in addition to condition \textup{(C4$'$)}, that
for each $y\in\bfS$ there exists $\bar{y}>y$ such that
$(y,\bar{y})\subseteq\bfS$ and
\[
\sup_{y\le y^\prime<\bar{y}} E_{x,y^\prime} \Bigl[ \sup_{t \geq N}e^{-qt}
\bigl|g(X_t)\bigr| \Bigr] \rightarrow0\qquad\mbox{as } N\uparrow\infty.
\]
Then, when $T=\infty$, $ v(x,\cdot)$ is right-continuous.
\end{proposition}

\begin{pf}
Choose $y\in\bfS$ and $y^\prime\in(y,\bar{y})$.
Applying Proposition \ref{equiv stop probl for diffu} with respect to
$x$ and $y^\prime$ yields
$v(x,y^\prime)=\sup_{\rho
\in\mathcal{M}} \tilde{E} [e^{-q\Gamma^\prime_\rho}
g(G_{\rho})]$ since $T=\infty$. Fix an arbitrary $\varepsilon>0$, and
choose an $\varepsilon$-optimal stopping time
$\rho^\prime_\varepsilon\in\mathcal{M}$ for $v(x,y^\prime)$ so that
%
%e4.7 #&#
\begin{equation}
\label{basic ineq} 0 \le v\bigl(x,y^\prime\bigr)-v(x,y) \le\varepsilon+
\tilde{E} \bigl[e^{-q\Gamma^\prime_{ \rho^\prime_\varepsilon}}g(G_{\rho
^\prime
_\varepsilon}) -e^{-q\Gamma_{ \rho^\prime_\varepsilon}}g(G_{\rho^\prime_\varepsilon})
\bigr].
\end{equation}
Because
$\Gamma_{ \rho^\prime_\varepsilon}\ge\Gamma^\prime_{ \rho^\prime
_\varepsilon}$,
the right-hand side of (\ref{basic ineq}) can be dominated by
\begin{eqnarray*}
&&\varepsilon+ \tilde{E} \bigl(1-e^{-q
(\Gamma_{ \rho^\prime_\varepsilon}-\Gamma^\prime_{ \rho^\prime
_\varepsilon})
} \bigr) e^{-q\Gamma^\prime_{ \rho^\prime_\varepsilon}}\bigl|g(G_{\rho^\prime
_\varepsilon})\bigr|
I\bigl(\rho^\prime_\varepsilon\le A^\prime_N
\bigr)
\\
&&\quad{}+ \tilde{E} e^{-q\Gamma^\prime_{ \rho^\prime_\varepsilon}}\bigl|g(G_{\rho^\prime
_\varepsilon})\bigr| I\bigl(
\rho^\prime_\varepsilon>A^\prime_N\bigr)
\\
&&\qquad\le\varepsilon+ \tilde{E} \bigl(1-e^{-q
(\Gamma_{ \rho^\prime_\varepsilon}-\Gamma^\prime_{ \rho^\prime
_\varepsilon})
} \bigr)
e^{-q\Gamma^\prime_{ \rho^\prime_\varepsilon}}\bigl|g(G_{\rho^\prime
_\varepsilon})\bigr| I\bigl(\rho^\prime_\varepsilon
\le A^\prime_N\bigr) \\
&&\qquad\quad{}+ \tilde{E} \Bigl[ \sup
_{t \geq N}e^{-qt} \bigl|g\bigl(\tilde{X}^\prime_t
\bigr)\bigr| \Bigr],
\end{eqnarray*}
where
\[
\tilde{E} \Bigl[ \sup_{t \geq N}e^{-qt}\bigl |g\bigl(
\tilde{X}^\prime_t\bigr)\bigl| \Bigr] = E_{x,y^\prime} \Bigl[
\sup_{t \geq N}e^{-qt} \bigl|g(X_t)\bigr| \Bigr]
\]
and both
%
%e4.8 #&#
\begin{eqnarray}
\label{lebesgue esti} e^{-q\Gamma^\prime_{
\rho^\prime_\varepsilon}}\bigl|g(G_{\rho^\prime
_\varepsilon})\bigr| &\le&\sup
_{t\le A_N^\prime} e^{-q\Gamma^\prime_{t}}\bigl|g(G_{t})\bigr|
\nonumber
\\
&\le&\sup_{t\le A_N^{\bar{y}}} e^{-q\Gamma^{\bar{y}}_{t}}\bigl|g(G_{t})\bigr|
\\
&=&\sup_{t\le A_N^{\bar{y}}} e^{-q\Gamma^{\bar{y}}_{t}}\bigl|
g\bigl(\tilde{X}^{\bar{y}}_{\Gamma^{\bar{y}}_t}
\bigr)\bigr| \le \sup_{t\le N}e^{-q{t}}\bigl|g\bigl(
\tilde{X}^{\bar{y}}_{t}\bigr)\bigr|\nonumber
\end{eqnarray}
and
%
%e4.9 #&#
\begin{eqnarray}
\label{conti esti} \Gamma_{ \rho^\prime_\varepsilon}-\Gamma^\prime_{ \rho^\prime
_\varepsilon} &=&
\int_0^{\rho^\prime_\varepsilon} \bigl(\xi_u^{-2}-
\bigl(\xi_u^\prime\bigr)^{-2}\bigr)\,du \le \int
_0^{A_N^\prime} \bigl(\xi_u^{-2}-
\bigl(\xi_u^\prime\bigr)^{-2}\bigr)\,du
\nonumber
\\[-8pt]
\\[-8pt]
\nonumber
& =&
\Gamma_{A_N^\prime}-N
\end{eqnarray}
on $\{\rho^\prime_\varepsilon\le A^\prime_N\}$.
Hence choosing $N$ large enough that
\[
\sup_{y\le y^\prime<\bar{y}} E_{x,y^\prime} \Bigl[ \sup_{t \geq N}e^{-qt}
\bigl|g(X_t)\bigr| \Bigr]\le\varepsilon,
\]
we obtain from (\ref{basic ineq})
%
%e4.10 #&#
\begin{equation}
\label{basic ineq2} \bigl|v\bigl(x,y^\prime\bigr)-v(x,y)\bigr| \le2\varepsilon+
\tilde{E} \bigl(1-e^{-q
(\Gamma_{ A_N^\prime}-N)} \bigr) \sup_{t\le N}e^{-q{t}}\bigl|g
\bigl(\tilde{X}^{\bar{y}}_{t}\bigr)\bigr|
\end{equation}
for some $N$ depending on $y$ but NOT on $y^\prime$.

Now, in inequality (\ref{basic ineq2}), replace $y$ and $y^\prime$
by $y_0$ and $y_n$, respectively,
with the $(y_n)$ bounded above by $\bar{y}$ and decreasing to $y_0$.
Since $\sup_{t\le N}e^{-q{t}}|g(\tilde{X}^{\bar{y}_0}_{t}) |$ is integrable,
it follows by dominated convergence that
\[
\lim_{n\to\infty}\bigl|v(x,y_n)-v(x,y_0)\bigr| \le2
\varepsilon+ \tilde{E} \bigl(1-e^{-q
(\lim_{n}\Gamma^0_{ A_N^n}-N)} \bigr) \sup_{t\le N}e^{-q{t}}\bigl|g
\bigl(\tilde{X}^{\bar{y}_0}_{t}\bigr)\bigr|,
\]
and so
\[
\lim_{n\to\infty}\bigl|v(x,y_n)-v(x,y_0)\bigr| \le2
\varepsilon
\]
by Lemma \ref{with feller}. Since $\varepsilon$ is arbitrary we conclude
with the desired result.
\end{pf}
%
%pr4.5 #&#
\begin{proposition}\label{conti finite time}
Assume, in addition to condition \textup{(C4$'$)}, that
the gain function $g$ is continuous.
Then, when $T<\infty$, $v(x,\cdot)$ is continuous.
\end{proposition}

\begin{pf}
Following the proof of the previous proposition
choose $y,y^\prime\in\bfS$ with $y<y^\prime$, fix an arbitrary
$\varepsilon
>0$ and
choose an $\varepsilon$-optimal stopping time
$\rho^\prime_\varepsilon\in\mathcal{M}^\prime_T$ so that
%
%e4.11 #&#
\begin{eqnarray}
\label{basic ineq new} 0 &\le& v\bigl(x,y^\prime\bigr)-v(x,y)
\nonumber
\\[-8pt]
\\[-8pt]
\nonumber
& \le&\varepsilon+
\tilde{E} \bigl[e^{-q\Gamma^\prime_{ \rho^\prime_\varepsilon}}g(G_{\rho
^\prime
_\varepsilon}) -e^{-q\Gamma_{ \rho^\prime_\varepsilon\wedge A_T}}g(G_{\rho^\prime
_\varepsilon\wedge A_T})
\bigr].
\end{eqnarray}
Note that $\rho^\prime_\varepsilon\le A^\prime_T$ and that
$\rho^\prime_\varepsilon\wedge A_T$ is used since one cannot conclude that
$v(x,y) \ge\tilde{E}e^{-q\Gamma_{ \rho}}g(G_{\rho})$
for stopping times $\rho$ which may exceed $A_T$ with
positive probability. Therefore, in contrast to the case where
$T=\infty$,
dominating the right-hand side of (\ref{basic ineq new}) leads to an
upper bound of
%
%e4.12 #&#
%e4.13 #&#
%e4.14 #&#
\begin{eqnarray}
&&\varepsilon+ \tilde{E} \bigl(1-e^{-q
(\Gamma_{ \rho^\prime_\varepsilon}-\Gamma^\prime_{ \rho^\prime
_\varepsilon})
} \bigr)
e^{-q\Gamma^\prime_{ \rho^\prime_\varepsilon}}\bigl|g(G_{\rho^\prime
_\varepsilon})\bigr| \label{eins}
\\
&&\qquad{}+ \tilde{E} \bigl(1-e^{-q
(T-\Gamma^\prime_{ \rho^\prime_\varepsilon})
} \bigr) e^{-q\Gamma^\prime_{ \rho^\prime_\varepsilon}}\bigl|g(G_{\rho^\prime
_\varepsilon})\bigr|
I\bigl(A_T<\rho^\prime_\varepsilon\le
A^\prime_T\bigr) \label{zwei}
\\
\label{drei} &&\qquad{}+ \tilde{E} e^{-qT}\bigl|g(G_{\rho^\prime_\varepsilon})-g(G_{A_T})\bigr| I
\bigl(A_T<\rho^\prime_\varepsilon\le A^\prime_T
\bigr)
\end{eqnarray}
by adding
$-e^{-qT}g(G_{\rho^\prime_\varepsilon})+e^{-qT}g(G_{\rho^\prime
_\varepsilon})$
in the case where $A_T<\rho^\prime_\varepsilon\le A^\prime_T$.

Now replace $y$ and $y^\prime$ by $y_n$ and $y_0$, respectively,
with $y_n\uparrow y_0$ in $\bfS$.
Suppose for now that Lebesgue's dominated convergence
theorem can be applied to interchange limit and expectation in
(\ref{eins}), (\ref{zwei}), (\ref{drei}).
Then it can be shown that
\[
\lim_{n\to\infty}\bigl|v(x,y_n)-v(x,y_0)\bigr| \le
\varepsilon
\]
proving left-continuity since $\varepsilon$ was arbitrary.
To see this first dominate
\[
\Gamma^n_{ \rho^0_\varepsilon}-\Gamma^0_{ \rho^0_\varepsilon}\quad
\mbox{by}\quad \Gamma^n_{A^0_T}-T
\]
performing a calculation similar to (\ref{conti esti}), but using $T$
instead of $N$. Then (\ref{eins}) tends to $\varepsilon$ as $n\to\infty$
by Lemma \ref{with feller}. Second, since
$\{A^n_T<\rho^0_\varepsilon\le A^0_T\}
=\{\Gamma^0_{A^n_T}<\Gamma^0_{\rho^0_\varepsilon}\le T\}$,
both (\ref{zwei}) and (\ref{drei}) converge to zero as $n\to\infty$
by Lemma \ref{with feller} and the continuity of $g$.

Finally it remains to justify the application of the dominated
convergence theorem. Observe that
\[
e^{-qT}\bigl|g(G_{\rho^0_\varepsilon})\bigr|
 \le e^{-q\Gamma^0_{ \rho^0_\varepsilon}}\bigl|g(G_{\rho^0_\varepsilon})\bigr|
= e^{-q\Gamma^0_{ \rho^0_\varepsilon}}\bigl|g\bigl(\tilde{X}^0_{\Gamma^0_{\rho
^0_\varepsilon}}\bigr)\bigr| \le \sup
_{0\leq t \leq T}e^{-qt} \bigl|g\bigl(\tilde{X}^0_t
\bigr)\bigr|
\]
since $\rho^0_\varepsilon\le A^0_T$ and
\[
e^{-qT}\bigl|g(G_{A^n_T})\bigr| = e^{-q\Gamma^n_{ A^n_T}}\bigl|g(G_{A^n_T})\bigr|
\le \sup_{t\le A_T^0}e^{-q\Gamma^0_{t}}\bigl|g(G_{t})\bigr| \le \sup
_{0\leq t \leq T}e^{-qt} \bigl|g\bigl(\tilde{X}^0_t
\bigr)\bigr|
\]
since $A^n_T\le A^0_T$ for all $n\ge1$
which, by (\ref{eq:integrabilityCond}),
gives an integrable bound
with respect to all three terms (\ref{eins}), (\ref{zwei}), (\ref{drei}).

For the right-continuity, replace $y$ and $y^\prime$ by $y_0$ and
$y_n$, respectively, assuming $y_n\downarrow y_0$ in $\bfS$. Note that
\[
e^{-qT}\bigl|g(G_{\rho^n_\varepsilon})\bigr| \le
e^{-q\Gamma^n_{ \rho^n_\varepsilon}}\bigl|g(G_{\rho^n_\varepsilon})\bigr|
\le \sup_{t\le T}e^{-q{t}}\bigl|g\bigl(\tilde{X}^1_{t}
\bigr)\bigr|,
\]
where the second inequality is obtained
following the line of inequalities in~(\ref{lebesgue esti}) but using
$T$ and $y_1$ instead of $N$ and $\bar{y}$, respectively.
As $e^{-qT}|g(G_{A^0_T})|\le\sup_{0\leq t \leq T}e^{-qt}|g(\tilde{X}^0_t)|$
too, dominated convergence can be applied again by (\ref{eq:integrabilityCond})
with respect to all three terms (\ref{eins}), (\ref{zwei}), (\ref{drei}).
Then (\ref{eins}) tends to $\varepsilon$ as $n\to\infty$ by Lemma \ref
{with feller}
since $\Gamma^0_{ \rho^n_\varepsilon}-\Gamma^n_{ \rho^n_\varepsilon}$
can be estimated by $\Gamma^0_{ A^n_T}-T$. Furthermore,
(\ref{zwei}) and (\ref{drei}) converge to zero as $n\to\infty$
by Lemma \ref{with feller} and the continuity of $g$ since
$T-\Gamma^n_{ \rho^n_\varepsilon}\le T-\Gamma^n_{ A^0_T}$ on
$\{A^0_T<\rho^n_\varepsilon\le A^n_T\}$.
So, making $\varepsilon$ arbitrarily small completes the proof.
\end{pf}

%s5 #&#
\section{Application to option pricing}\label{option pricing}

Assume that the dynamics of $X$ are given by
\renewcommand{\theequation}{1.\arabic{equation}$'$}
\setcounter{equation}{0}
\begin{equation}\label{eq1.1prime}
dX=XY\,dB,%\leqno(1.1^\prime)
\end{equation}
which is the special case $a(x)=x$ of equation (\ref{eq:dynX}).
In mathematical finance~(\ref{eq1.1prime}) describes a simple model for the
discounted price of an asset with stochastic volatility $Y$.

If exercised at a stopping time $\tau$, the American options we have
in mind would pay off $g(e^{r\tau}X_\tau)$ where $r>0$ stands for the
instantaneous interest rate which is assumed to be constant. So, for
notational convenience, the discount rate $q$ is replaced by $r$
throughout this section.

In this setup, assuming the measure $P_{x,y}$ is used for pricing
when $X_0=x$ and $Y_0=y$,
the price of such an option with maturity $T\in[0,\infty]$ is
\begin{equation}%\label{eq1.2prime}
v(x,y)=\sup_{0\leq\tau\leq T} E_{x,y} \bigl[ e^{- r\tau}g
\bigl(e^{r\tau
}X_\tau\bigr)\bigr], %\leqno(1.2^\prime)
\end{equation}
where the supremum is taken over all finite stopping times with
respect to the filtration generated by $(X,Y)$.
This value function differs from the value function given
by (\ref{osp}) since $g$ is not applied to $X_\tau$ but to $e^{r\tau
}X_\tau$
and, as a consequence, some of the conditions for our results have to
be adjusted slightly.

First, the condition
\begin{equation}\label{eq1.3prime}\quad
E_{x,y} \Bigl[ \sup_{0\leq t \leq T}e^{-rt} \bigl|g
\bigl(e^{rt}X_t\bigr)\bigr| I(t<\infty) \Bigr]< \infty \qquad\mbox{for
all } (x,y)\in\bfR\times\bfS %\leqno(1.3^\prime)
\end{equation}
is now assumed throughout. Then
\[
v(x,y)=\sup_{\tau\in\mathcal{T}_T } \tilde{E} \bigl[e^{-r \tau}
g\bigl(e^{r\tau} \tilde{X}_\tau\bigr)\bigr] =\sup
_{\rho\in\mathcal{M}_T} \tilde{E} \bigl[e^{-r\Gamma_\rho}g\bigl(e^{r\Gamma_\rho}G_{\rho}
\bigr)\bigr]
\]
is the analogue to what was obtained in
Propositions \ref{equiv stop probl} and \ref{equiv stop probl for diffu}
for the value function given by (\ref{osp}).
However, in order to conclude the results
of Theorems \ref{main one} and \ref{main two}
for the new value function,
a new condition has to be imposed on $g$.
%
%co5.1 #&#
\begin{corollary} \label{cor:mathfinapp}
Let $v$ be the value function given by (\ref{eq1.2prime}).
In addition to the assumptions made in
either Theorem \ref{main one} or \ref{main two}
assume that $g$ is a decreasing function.
Define ${\cal K}^{g+}_T$ to be the collection
of all finite stopping times $\tau\le T$
with respect to the filtration generated by $(X,Y)$
such that $g(e^{r\tau}X_\tau)\ge0$.
Fix $(x,y)\in\bfR\times\bfS$ and assume that
$v(x,y) = \sup_{\tau\in{\cal K}^{g+}_T}E_{x,y} [
e^{-r\tau}g(e^{r\tau}X_\tau)]$.
Then
\[
v(x,y)\le v\bigl(x,y^\prime\bigr) \qquad\mbox{for all $y^\prime\in
\bfS$ such that $y\le y^\prime$,}
\]
so that $v(x,y)$ is a lower bound for $v(x,\cdot)$ on
$[y,\infty)\cap\bfS$.
\end{corollary}
%
%re5.2 #&#
\begin{remark}\label{new g}
(i) The proofs of this and the next corollary are contained in
the \hyperref[app]{Appendix}.\vspace*{-6pt}
\begin{longlist}[(iii)]
\item[(ii)] If $g$ is a monotone function, then it has a left and a
right-continuous
version. Note that the proof of Corollary \ref{cor:mathfinapp}
does not depend on choosing a specific version for $g$.
But, when applying the corollary to show continuity properties of the value
function, we will choose the right-continuous version in what follows.

\item[(iii)] Of course, Corollary \ref{cor:mathfinapp} does not depend
on the specific
choice of the diffusion coefficient $a$ in this section
as long as (\ref{eq1.3prime}) and all other assumptions of
Theorems \ref{main one} or \ref{main two} are satisfied.

\item[(iv)] If $a(x)=x$, then
conditions (C2) or (C2$'$)
assumed in Corollary \ref{cor:mathfinapp} ensures that the
discounted price $X$ is a positive exponential local martingale
of the form
\[
X_t=x\exp\biggl\{\int_0^t
Y_s \,dB_s - \frac{1}{2}\int_0^t
Y_s^2 \,ds\biggr\},\qquad t\ge0, P_{x,y}\mbox{-a.s.},
\]
since the stochastic integrals $\int_0^t Y_s \,dB_s, t\ge0$,
are all well defined.
Furthermore, because
$\lim_{t\uparrow\infty} \int_0^t Y_s^2 \,ds = \infty$ $P_{x,y}$-a.s.,
$X_t$ tends to zero for large $t$ as in the Black--Scholes model.

\item[(v)] From (iv) above it follows immediately that, in the case $a(x)=x$,
all processes satisfying conditions (C1) and (C2) on page
\pageref{a and b} have the
same law.

\item[(vi)] Note that, in this section,
the equation for $G$ in (\ref{equ star}) on page \pageref{equ star}
coincides with the linear equation $dG=G\,dW$ which has a unique
nonexploding strong solution for all $G_0\in\bfR$. Hence
condition (C3$'$) on page \pageref{only c} becomes a condition only
on the
coefficients $\eta,\theta$ of the equation for $\xi$ in (\ref{equ star}).
\end{longlist}
\end{remark}

We now consider the diffusion case and
discuss the results of Section~\ref{conti}
for the value function given by (\ref{eq1.2prime}).
So, let $\bfS$ be an open subset of $(0,\infty)$, fix $x\in\bfR$ and
replace condition (C4$'$) on page \pageref{summarising} by:
\begin{itemize}[(C4$''$):]
\item[(C4$''$):]
\begin{itemize}
\item
the gain function $g$ is decreasing and
satisfies $\{g\ge0\}\neq\varnothing$;
\item
the process $(X,Y)$ satisfies conditions
(\ref{eq1.3prime}), (C1$'$), (C2$'$), and the value
function $v$ satisfies
\renewcommand{\theequation}{4.1$'$}
\begin{equation}\label{eq4.1prime}
v(x,y) = \sup_{\tau\in{\cal
K}^{g+}_T}E_{x,y} \bigl[e^{-r\tau}g
\bigl(e^{r\tau}X_\tau\bigr)\bigr] \qquad\mbox{for all $y\in\bfS$}
\end{equation}
for the chosen $x$
(using the definition of ${\cal K}^{g+}_T$
given in Corollary~\ref{cor:mathfinapp});\vadjust{\goodbreak}
\item
condition (C3$'$) holds true for system (\ref{equ star}), and the
second equation in~(\ref{equ star}) has a Feller solution.
\end{itemize}
\end{itemize}
%
%co5.3 #&#
\begin{corollary} \label{rem:contPcor}
Let $v$ be the value function given by (\ref{eq1.2prime}).
Assume condition \textup{(C4$''$)}.
\begin{longlist}[(iii)]
\item[(i)] If $g$ is bounded from below then, when $T=\infty$,
$v(x,\cdot)$ is left-continuous and lower semicontinuous.

\item[(ii)]
If $g$ is continuous and if
for each $y\in\bfS$ there exists $\bar{y}>y$ such that
$(y,\bar{y})\subseteq\bfS$ and
\[
\sup_{y\le y^\prime<\bar{y}} E_{x,y^\prime} \Bigl[ \sup_{t \geq N}e^{-rt}
\bigl|g\bigl(e^{rt}X_t\bigr)\bigr| \Bigr] \rightarrow0\qquad\mbox{as }N
\uparrow\infty,
\]
then, when $T=\infty$, $v(x,\cdot)$ is right-continuous.

\item[(iii)] If $g$ is bounded from below and continuous
then, when $T<\infty$, $v(x,\cdot)$ is continuous.
\end{longlist}
\end{corollary}

%s6 #&#
\section{Examples}\label{example section}

We now discuss three models used in option pricing and explain the
impact of our results.

\textit{Pricing of American Puts via Jobert and Rogers} \cite{J-R}
using the Markov modulated model
\[
dX=XY\,dB,\qquad\mbox{$Y$ finite state Markov chain}.
\]
Notice that the value function in \cite{J-R} is
more general than ours as the authors allow for an interest rate which
depends on $Y$. So in what follows we always mean a constant interest rate
when applying our results to the value function\footnote{Note that
the notation of the value function in \cite{J-R} is different because
our Markov
chain $Y$ is, in their terms, a function $\sigma$ applied to the
Markov chain playing the role of their volatility process.}
in \cite{J-R}.

Obviously, the gain function
$g(x)=\max\{0,K-x\}$ where $K$ is the strike price
is decreasing and satisfies both condition (\ref{eq1.3prime}) and
\[
v(x,y) = \sup_{\tau\in{\cal K}^{g+}_T}E_{x,y} \bigl[ e^{-r\tau}g
\bigl(e^{r\tau}X_\tau\bigr)\bigr],\qquad (x,y)\in\bfR\times\bfS.
\]
So, recalling Remark \ref{new g}(iv)${}+{}$(v),
Corollary \ref{cor:mathfinapp} implies that, for fixed $x\in\bfR$, the
value function
$v(x,y)$ in \cite{J-R} is monotonously increasing in $y\in\bfS=\{
y_1,\ldots,y_m\}$,
provided $Y$ is skip-free.

Knowing this monotonicity property of the value function massively reduces
the computational complexity of PROBLEM 1 on page 2066 in \cite{J-R}.
The authors verified that the value function is uniquely attained
at a stopping time of the form\footnote{We again adapted the author's
notation to ours in the definition of $\tau^\star$.}
\[
\tau^\star= \inf\bigl\{t\ge0\dvtx X_t<b[Y_t]
\bigr\},
\]
where the vector $b[y_i], i=1,\ldots,m$, is indexed by the states of
the Markov chain $Y$ and their PROBLEM 1 consists in finding the
so-called \textit{thresholds} $b[y_i]$ which are assumed to be in the
order $b[y_1]\ge\cdots\ge b[y_m]$.
It is then stated in a footnote on the same page, 2066, that
``When it comes in practice to identifying the thresholds,
no assumption is made on the ordering, and all possible orderings are
considered.''
Of course, this approach has exponential complexity. Our result on the
monotonicity of the value function would reduce this complexity to choosing
one ordering $b[y_1]>\cdots>b[y_m]$ if $y_1<\cdots<y_m$ and
$Y$ is skip-free. Indeed, since $\tau^\star$ is the
unique optimal stopping time for this problem, by general theory,
it must coincide with the first time the process $(X,Y)$ enters the
stopping region $\{(x,y)\dvtx v(x,y)=g(x)\}$.
Thus,
as it is not optimal to stop when $g$ is zero, we obtain that
\[
\mbox{$v(x,y_i)=g(x)$ for $x\le b[y_i]$} \quad\mbox{while}\quad
\mbox{$v(x,y_i)>g(x)$ for $x>b[y_i]$}
\]
for each $i=1,\ldots,m$ which gives the unique ordering of the
thresholds since $g$ is strictly decreasing on $\{g>0\}$.

\textit{The Hull and White model} \cite{HW}:
\[
dX = X \sqrt{V} \,dB \quad\mbox{and}\quad \,dV = 2\eta V \,dB^Y + \kappa V \,dt,
\]
where $\eta,\kappa>0$
and $B,B^Y$ are independent Brownian motions.\footnote{Remark that
$\langle B,B^Y\rangle\neq0$ is possible but we follow
Hull and White's original setup.}
Setting $Y=\sqrt{V}$ transforms the above system into
\[
dX=XY\,dB \quad\mbox{and}\quad dY = \eta Y \,dB^Y + \theta Y \,dt,
\]
where $\theta=(\kappa-\eta^2)/2$.
Assuming a positive initial condition, this equation has a
pathwise unique positive solution for every $\eta,\theta\in\bfR$.
Calculating the equation for $\xi$ in (\ref{equ star}) on page
\pageref
{equ star}
gives a constant diffusion coefficient $\eta$, and
if $Z$ denotes $\xi/\eta$, then
\[
dZ=dW^\xi+\frac{\theta}{\eta^2}Z^{-1}\,dt,
\]
which formally is an equation
for a Bessel process of dimension $\phi=1+2\theta/\eta^2$.
This equation, and so the equation for $\xi$, only has a unique
nonexploding strong solution if $\phi\geq2$, and this
solution stays positive when started from a positive initial
condition. As made clear in Section~\ref{diffu},
the fact that $Y$ satisfies
condition~(C2$'$) on page \pageref{a,b and c} can be derived
from condition (\ref{inftyprop}) with respect to
\[
\Gamma_t=\int_0^t
\frac{1}{\xi_u^2}\,du=\eta^2\int_0^t
\frac{1}{Z_u^2}\,du, \qquad t\ge0.
\]
Now, by applying Proposition A.1(ii)--(iii) in \cite{Hobson} with
respect to
the second time integral above, we see that $\Gamma$ satisfies
condition (\ref{inftyprop}) if $\phi\geq2$.
So, assuming $\phi\geq2$, Remark \ref{new g}(iv)${}+{}$(vi) ensures that
there is a unique strong Markov process $(X,Y)$
which satisfies conditions (C1$'$) and (C2$'$) on page
\pageref{a,b and c} and that the system (\ref{equ star}) satisfies
condition (C3$'$) on page \pageref{only c} in this example.
Since Bessel processes are Feller processes (see \cite{Rev-Yor}, page 446),
the second equation of (\ref{equ star})
has a Feller solution.

Therefore if $\phi\geq2$ (i.e., $\kappa\geq2\eta^2$), then the
conclusions of Corollaries \ref{cor:mathfinapp} and~\ref{rem:contPcor}
apply to perpetual American options whenever the
corresponding pay-off function $g$ satisfies the conditions stated.

\textit{The Heston model} \cite{Heston}:
\[
dX = X\sqrt{V} \,dB \quad\mbox{and} \quad dV = 2\eta\sqrt{V} \,dB^Y + \kappa(
\lambda-V) \,dt,
\]
where $\eta,\kappa,\lambda>0$ are constants,
and $B,B^Y$ are Brownian motions, this time with covariation
$\delta\in[-1,1]$.
The equation for $V$ describes the so-called Cox--Ingersoll--Ross process,
and it is well known (see \cite{CIR}, page 391) that,
with a positive initial condition, this\vadjust{\goodbreak}
equation has a pathwise unique positive solution
if $\kappa\lambda\ge2\eta^2$.
Setting $Y=\sqrt{V}$ transforms the system into
\[
dX=XY\,dB\quad \mbox{and}\quad dY = \eta \,dB^Y + \biggl( \frac{\theta_1}{Y}-
\theta_2 Y \biggr)\,dt
\]
with $\theta_1=(\kappa\lambda-\eta^2)/2$ and $\theta_2=\kappa/2$.
It is clear that the pathwise uniqueness of the equation for $V$
ensures the pathwise uniqueness of positive solutions of the equation
for $Y$.
Calculating the equation for $\xi$ in (\ref{equ star}) on page
\pageref{equ star} yields
\[
d\xi= \frac{\eta}{\xi} \,dW^\xi+ \biggl(\frac{\theta_1}{\xi
^3}-
\frac
{\theta_2}{\xi} \biggr) \,dt,
\]
and
hence $Z=\xi^2/(2\eta)$ satisfies
\[
dZ = dW^\xi+ \biggl( \frac{\phi-1}{2 Z} - \frac{\theta_2}{\eta
}
\biggr)\,dt
\]
with $\phi=\theta_1/\eta^2 + 3/2$.
By changing to an equivalent probability measure,
this equation for $Z$ is transformed into an equation
for a Bessel process of dimension $\phi$ which only has a unique
nonexploding strong solution if $\phi\ge2$, and this
unique strong solution stays positive when started from a positive
initial condition.
All these properties and the Feller property of Bessel
processes carry over to the solutions of the equation for $\xi$.
Finally, the process
\[
\Gamma_t=\int_0^t
\frac{1}{\xi_u^2}\,du = \frac{1}{2\eta}\int_0^t
\frac{1}{Z_u}\,du, \qquad t\ge0,
\]
satisfies (\ref{inftyprop})
if $\phi\geq2$ (apply Proposition A.1(ii)--(iii) in \cite{Hobson}
to the second integral)
which implies condition (C2$'$) on page \pageref{a,b and c} following
the arguments
given in Section~\ref{diffu}.
So,
as in the previous example, all conditions imposed on $X,Y,\xi$ in
the Corollaries \ref{cor:mathfinapp} and \ref{rem:contPcor} are
satisfied if $\phi\geq2$ or equivalently $\kappa\lambda\geq2\eta^2$.

%s7 #&#
 \begin{appendix}\label{app}
\section*{Appendix}

\begin{pf*}{Proof of Lemma \ref{lem_equiv_filt}}
Fix $\rho\in\mathcal{M}$ and $r\ge0$, and set
\[
\Omega_0 = \bigl\{\omega\in\Omega: \mbox{$s<\Gamma_t(
\omega)$ if and only if $A_s(\omega)<t$ for all $0\leq s,t<\infty$}
\bigr\}.
\]
Then
\[
\{\Gamma_\rho\leq r\}\cap\Omega_0\cap\{A_r<
\infty\} = \{\rho\leq A_r\}\cap\Omega_0\cap
\{A_r<\infty\}
\]
implies
\[
\{\Gamma_\rho\leq r\}\in{\cal F}^{G,Z,Z^\prime}_{A_r}
\]
since both
$P(\Omega_0\cap\{A_r<\infty\})=1$ by property (P2), (\ref{inftyprop
for A})
and
$\{\rho\leq A_r\}\in{\cal F}^{G,Z,Z^\prime}_{A_r}$. Note that
$\Omega_0\cap\{A_r<\infty\}\in{\cal F}^{G,Z,Z^\prime}_{A_r}$ as
${\cal F}^{G,Z,Z^\prime}_0$ already contains all $\tilde{P}$-null
sets.\vadjust{\goodbreak}

Similarly, if $\tau\in\mathcal{T}$, then
$\{A_\tau\leq r\}\in{\cal F}^{G,Z,Z^\prime}_{A_{\Gamma_r}}$
where $A_{\Gamma_r}=r$ a.s. by property~(P1). Thus the inclusion
${\cal F}^{G,Z,Z^\prime}_{A_{\Gamma_r}}\subseteq{\cal
F}^{G,Z,Z^\prime}_r$
must be true.
\end{pf*}

\begin{pf*}{Proof of (\ref{with snell})}
\label{snelli}
By (\ref{one prob space}), we only have to show that
%
%e7.1 #&#
\renewcommand{\theequation}{\Alph{section}.\arabic{equation}}
\begin{equation}
\label{without theta} \sup_{ 0\le\tilde{\tau}\le T } \tilde{E} \bigl[e^{-q\tilde{\tau
}}g(
\tilde {X}_{\tilde{\tau}})\bigr] = \sup_{ 0\le\tau\le T } \tilde{E}
\bigl[e^{-q\tau}g(\tilde{X}_{\tau})\bigr],
\end{equation}
where $\tilde{\tau}$ on the above left-hand side corresponds to
finite stopping times with respect to the filtration ${\cal F}^{\tilde
{X},\tilde{Y}}_t, t\ge0$,
generated by the pair of processes $(\tilde{X},\tilde{Y})$ while
$\tau$ on the above right-hand side corresponds to
finite stopping times with respect to the possibly bigger filtration
${\cal F}_{A_t}^{G,Z,Z'}, t\ge0$.
In what follows we assume that ${\cal F}^{\tilde{X},\tilde{Y}}_t,
t\ge
0$, was augmented.
Without loss of generality,
we also assume that there exist a family $\{\theta_t, t\ge0\}$
of shift operators on our chosen probability space
$(\tilde{\Omega}, \tilde{\mathcal{F}},\tilde{P})$.

We are going to show that
\[
\sup_{ \tilde{\tau}\in{\cal O}_0^T} \tilde{E} \bigl[e^{-q\tilde{\tau}}g(
\tilde{X}_{\tilde{\tau}})\bigr] = \sup_{ 0\le\tau\le T } \tilde{E}
\bigl[e^{-q\tau}g(\tilde{X}_{\tau})\bigr],
\]
where ${\cal O}_s^T$ stands for the family of
all finite ${\cal F}^{\tilde{X},\tilde{Y}}_t$-stopping times
$\tilde
{\tau}$
satisfying $s\le\tilde{\tau}\le T$ and
$\tilde{\tau}-s\stackrel{\mathrm{a.s.}}{=}\gamma\circ\theta_s$
for some ${\cal F}^{\tilde{X},\tilde{Y}}_\infty$-measurable random
variable $\gamma$.
This obviously proves (\ref{without theta}) because the above left-hand
side is
less than or equal to the left-hand side of (\ref{without theta}).

First observe that
\[
{\cal F}^{\tilde{X},\tilde{Y}}_0 = {\cal F}_{0}^{G,Z,Z'}
= {\cal F}_{A_0}^{G,Z,Z'} = \sigma\qquad (\tilde{P}\mbox{-null sets})
\]
hence
\begin{eqnarray}
\sup_{ \tilde{\tau}\in{\cal O}_0^T} \tilde{E} \bigl[e^{-q\tilde{\tau}}g(
\tilde{X}_{\tilde{\tau}})\bigr] \stackrel{\mathrm{a.s.}} {=}
\tilde{V}_0 \nonumber\\
\eqntext{\mbox{where }\displaystyle \tilde{V}_t =
\esssup_{\tilde{\tau}\in{\cal O}_t^T} \tilde{E} \bigl[ e^{-q\tilde{\tau}}g(\tilde{X}_{\tilde{\tau}})
| {\cal F}^{\tilde
{X},\tilde{Y}}_t \bigr]}
\end{eqnarray}
and
\begin{eqnarray}
\sup_{ 0\le{\tau}\le T} \tilde{E} \bigl[e^{-q{\tau}}g(
\tilde{X}_{{\tau}})\bigr] \stackrel{\mathrm{a.s.}} {=} {V}_0\nonumber\\
\eqntext{\mbox{where } \displaystyle{V}_t = \esssup_{t\le{\tau}\le T} \tilde{E} \bigl[
e^{-q{\tau}}g(\tilde{X}_{{\tau}}) | {\cal F}_{A_t}^{G,Z,Z'}
\bigr].}
\end{eqnarray}
Note that $t\in{\cal O}_t^T$ gives
$\tilde{V}_t \ge e^{-q{t}}g(\tilde{X}_{{t}})$ almost surely
for each $t\ge0$.

Second, since
$(\tilde{X},\tilde{Y})$ has the same law as $({X},{Y})$ under $P_{x,y}$,
the process $(\tilde{X},\tilde{Y})$ is strong Markov with respect to
${\cal F}^{\tilde{X},\tilde{Y}}_t, t\ge0$. Therefore
%
%e7.2 #&#
\begin{eqnarray}
\label{marki} %
\tilde{E} \bigl[
e^{-q\tilde{\tau}}g(\tilde{X}_{\tilde{\tau}}) | {\cal F}^{\tilde
{X},\tilde{Y}}_t
\bigr] &\stackrel{\mathrm{a.s.}} {=}& \tilde{E} \bigl[ e^{-q\tilde{\tau}}g(
\tilde{X}_{\tilde{\tau}}) | \sigma(\tilde {X}_t,\tilde{Y}_t)
\bigr]
\nonumber
\\[-8pt]
\\[-8pt]
\nonumber
&\stackrel{\mathrm{a.s.}} {=}& \tilde{E} \bigl[
e^{-q\tilde{\tau}}g(\tilde{X}_{\tilde{\tau}}) | {\cal F}_{A_t}^{G,Z,Z'}
\bigr]
\end{eqnarray}
for all $\tilde{\tau}\in{\cal O}_t^T$ because $(\tilde{X},\tilde{Y})$
is strong Markov with respect to ${\cal F}_{A_t}^{G,Z,Z'}, t\ge0$, too.
Note that we have only used the Markov property to get (\ref{marki}).

Third, (\ref{marki}) implies that $\tilde{V}_t, t\ge0$, is an
${\cal F}_{A_t}^{G,Z,Z'}$-supermartingale.
Using assumption (\ref{eq:integrabilityCond}),
the proof of this fact is almost identical to part $1^o$ of the proof of
Theorem 2.2 in \cite{PS2006}.
The only difference is concerned with ${\cal F}_{A_t}^{G,Z,Z'}$-stopping times of type
\[
\tau= \tilde{\tau}_1 I_\Gamma+\tilde{\tau}_2
I_{\tilde{\Omega
}\setminus\Gamma}
\]
given by
\[
\Gamma= \bigl\{ \tilde{E} \bigl[ e^{-q\tilde{\tau}_1}g(\tilde{X}_{\tilde{\tau}_1}) |
{\cal F}_{A_t}^{G,Z,Z'} \bigr] \ge \tilde{E} \bigl[
e^{-q\tilde{\tau}_2}g(\tilde{X}_{\tilde{\tau}_2}) | {\cal F}_{A_t}^{G,Z,Z'}
\bigr] \bigr\}
\]
and $\tilde{\tau}_1,\tilde{\tau}_2\in{\cal O}_t^T$ where $t\ge0$
is fixed.
We need to show that $\tau\in{\cal O}_t^T$.
But, if
$\tilde{\tau}_1-t\stackrel{\mathrm{a.s.}}{=}\gamma_1\circ\theta_t$
and
$\tilde{\tau}_2-t\stackrel{\mathrm{a.s.}}{=}\gamma_2\circ\theta_t$,
then, by (\ref{marki}),
%
%e7.3 #&#
\begin{equation}
\label{shifti} \tau(\omega)-t = \gamma_1(\theta_t
\omega) I_{\Gamma_0}(\theta_t\omega) + \gamma_1(
\theta_t\omega) I_{\tilde{\Omega}\setminus\Gamma
_0}(\theta _t\omega)
\end{equation}
for almost every $\omega\in\tilde{\Omega}$. Here $\Gamma_0$ stands for
a set of type
$\{\phi_1(\tilde{X}_0,\tilde{Y}_0)\ge\phi_2(\tilde{X}_0,\tilde
{Y}_0)\}$ where
$\phi_1,\phi_2\dvtx \bfR^2\to\bfR$ are Borel-measurable functions satisfying
\[
\phi_1(\tilde{X}_t,\tilde{Y}_t) \stackrel{
\mathrm{a.s.}} {=} \tilde{E} \bigl[ e^{-q\tilde{\tau}_1}g(\tilde{X}_{\tilde{\tau}_1}) |
\sigma(\tilde {X}_t,\tilde{Y}_t) \bigr]
\]
and
\[
\phi_2(\tilde{X}_t,\tilde{Y}_t) \stackrel{
\mathrm{a.s.}} {=} \tilde{E} \bigl[ e^{-q\tilde{\tau}_2}g(\tilde{X}_{\tilde{\tau}_2}) |
\sigma(\tilde {X}_t,\tilde{Y}_t) \bigr],
\]
and hence (\ref{shifti}) justifies $\tau\in{\cal O}_t^T$.

Now, by Theorem 2.2 in \cite{PS2006}, the Snell envelope $V_t, t\ge0$,
is the smallest ${\cal F}_{A_t}^{G,Z,Z'}$-supermartingale
dominating the gain process and hence $\tilde{V}_t\ge V_t$ almost surely
for each $t\ge0$
proving
\[
\sup_{ \tilde{\tau}\in{\cal O}_0^T } \tilde{E} \bigl[e^{-q\tilde{\tau
}}g(
\tilde{X}_{\tilde{\tau}})\bigr] \ge \sup_{ 0\le\tau\le T } \tilde{E}
\bigl[e^{-q\tau}g(\tilde{X}_{\tau})\bigr].
\]
The reverse inequality is obvious.
\end{pf*}

\begin{lemmas}\label{A1}
Let $W,Z,A,\tilde{M}$ be given on the filtered probability space
$(\tilde{\Omega}, \tilde{\mathcal{F}},\tilde{\cal F}_t,t\ge
0,\tilde{P})$
as introduced in Section~\ref{regswi}. Then the time-changed processes
$\tilde{B}=\tilde{M}\circ A$ and $\tilde{Y}=Z\circ A$ are independent.
\end{lemmas}
\begin{pf}
Let $\mathcal{F}_t^W, t\geq0$, denote the augmentation
of the filtration generated by $W$ and
define the so-called big filtration by
\[
{\mathcal{F}}_t^{\mathrm{big}} = \mathcal{F}_t^W
\vee \sigma\bigl(\{Z_s\dvtx s\ge0\}\bigr), \qquad t\ge0.\vadjust{\goodbreak}
\]
Note that $W$ is an ${\mathcal{F}}^{\mathrm{big}}_t$ Brownian motion
since $W$ and $Z$ are independent, and hence the stochastic integral
$\tilde{M}$ is a continuous ${\mathcal{F}}^{\mathrm{big}}_t$ local martingale.
Since $A$ is a functional of $Z$, it must be an
${\mathcal{F}}^{\mathrm{big}}_t$ time-change by the definition of the
big filtration. As $A$ satisfies (\ref{inftyprop for A}), it follows from
Dambis--Dubins--Schwarz' theorem~\cite{Rev-Yor}, Theorem~V.1.6, that
$\tilde{B}=\tilde{M}\circ A$ is an ${\mathcal{F}}^{\mathrm{big}}_{A_t}$
Brownian motion.
But $\tilde{Y}=Z\circ A$ is a functional of $Z$, so it must be
independent of $\tilde{B}$ since
$\sigma(\{Z_s\dvtx s\ge0\})\subseteq
{\mathcal{F}}^{\mathrm{big}}_0=\tilde{\mathcal{F}}^{\mathrm{big}}_{A_0}$,
and $\tilde{B}$ is independent of ${\mathcal{F}}^{\mathrm{big}}_{A_0}$.
\end{pf}

\begin{pf*}{Proof of Corollary \ref{cor:mathfinapp}}
The only part of the proof
where the additional condition on $g$ is needed
is the verification of (\ref{compari}).
But, for (\ref{eq1.2prime}),
the modification of (\ref{compari}) reads
\[
\tilde{E} [e^{-r\Gamma_\rho}g\bigl(e^{r\Gamma_\rho} G_\rho\bigr) \le
\tilde{E} \bigl[e^{-r\Gamma^\prime_\rho}g\bigl(e^{r\Gamma^\prime_\rho} G_\rho\bigr)
\bigr]\qquad \mbox{for every $\rho\in\mathcal{M}_T^+$},
\]
and the above inequality is indeed true
because $\Gamma_t\ge\Gamma^\prime_t, t\ge0$, a.s.,
and $g$ is decreasing.
Note that the above set of stopping times $\mathcal{M}_T^+$
now denotes the set
$\{\rho\in\mathcal{M}_T\dvtx g(e^{r\Gamma_\rho} G_\rho)\ge0 \mbox{ a.s.}\}$.
\end{pf*}
\begin{pf*}{Proof of Corollary \ref{rem:contPcor}}
First observe that Lemma \ref{with feller} follows by simply
applying Corollary \ref{cor:mathfinapp} instead of Theorem \ref{main two}
and can therefore be used in the proof below.

Now, as the left-hand side of the estimate (\ref{ineq})
is trivially bounded from below since $g$ is bounded from below,
we obtain
\[
\tilde{E} e^{-r\Gamma^0_\rho} g\bigl(e^{r\Gamma^0_\rho}G_\rho\bigr) \leq
\liminf_{n\rightarrow\infty} \tilde{E} e^{-r\Gamma^n_\rho}g\bigl(e^{r\Gamma^n_\rho}G_{\rho}
\bigr)
\]
using Fatou's lemma, Lemma \ref{with feller} and Remark \ref{new g}(ii).
The remaining arguments below (\ref{fatou}) used to show
Proposition \ref{left conti for inf} also apply in the case
where (\ref{eq1.2prime}) holds proving the left-continuity claimed in part (i).
And finally,
the lower semicontinuity follows by the argument for lower semicontinuity
given in Remark \ref{lower semi cont}(i).

The proof of part (ii) is along the lines of the proof of
Proposition \ref{new int cond} with some small changes
emphasized below.

First, using the value function defined in (\ref{eq1.2prime}),
the right-hand side of (\ref{basic ineq}) is dominated by
\begin{eqnarray*}
&&\varepsilon+ \tilde{E} \bigl(1-e^{-r
(\Gamma_{ \rho^\prime_\varepsilon}-\Gamma^\prime_{ \rho^\prime
_\varepsilon})
} \bigr) e^{-r\Gamma^\prime_{ \rho^\prime_\varepsilon}}\bigl |g
\bigl(e^{r\Gamma^\prime_{ \rho^\prime_\varepsilon}}G_{\rho^\prime
_\varepsilon}\bigr)\bigr| I\bigl(\rho^\prime_\varepsilon
\le A^\prime_N\bigr)
\\
&&\qquad{}+ \tilde{E} e^{-r\Gamma_{ \rho^\prime_\varepsilon}}
\bigl|g\bigl(e^{r\Gamma^\prime_{ \rho^\prime_\varepsilon}}G_{\rho^\prime
_\varepsilon}\bigr)
-g\bigl(e^{r\Gamma_{ \rho^\prime_\varepsilon}}G_{\rho^\prime_\varepsilon}\bigr)\bigr| I\bigl(\rho^\prime_\varepsilon
\le A^\prime_N\bigr)
\\
&&\qquad{} + \tilde{E} \Bigl[ \sup_{t \geq N}e^{-rt}
\bigl|g\bigl(\tilde{X}^\prime_t\bigr)\bigr| \Bigr] + \tilde{E} \Bigl[
\sup_{t \geq N}e^{-rt}\bigl |g(\tilde{X}_t)\bigr| \Bigr],
\end{eqnarray*}
where the middle term
%
%e7.4 #&#
\renewcommand{\theequation}{\Alph{section}.\arabic{equation}}
\begin{equation}
\label{extra term} \tilde{E} e^{-r\Gamma_{ \rho^\prime_\varepsilon}} \bigl|g\bigl(e^{r\Gamma^\prime_{ \rho^\prime_\varepsilon}}G_{\rho^\prime
_\varepsilon}
\bigr) -g\bigl(e^{r\Gamma_{ \rho^\prime_\varepsilon}}
G_{\rho^\prime_\varepsilon}\bigr)\bigr| I\bigl(
\rho^\prime_\varepsilon\le A^\prime_N\bigr)
\end{equation}
is new.
Note that the $\varepsilon$-optimal stopping time
$\rho^\prime_\varepsilon$ can be chosen from the set
$(\mathcal{M}^\prime)^+=\{\rho\in{\cal M}\dvtx g(e^{r\Gamma^\prime_{ \rho}}G_{\rho})\ge0\}$
and so
\begin{eqnarray*}
e^{-r\Gamma^\prime_{ \rho^\prime_\varepsilon}}
\bigl|g\bigl(e^{r\Gamma^\prime_{ \rho^\prime_\varepsilon}}G_{\rho^\prime
_\varepsilon}\bigr)\bigr| &=&
e^{-r\Gamma^\prime_{ \rho^\prime_\varepsilon}} g\bigl(e^{r\Gamma^\prime_{ \rho^\prime_\varepsilon}}G_{\rho^\prime
_\varepsilon}\bigr)
\\
&\le& e^{-r\Gamma^{\bar{y}}_{ \rho^\prime_\varepsilon}} g\bigl(e^{r\Gamma^{\bar{y}}_{ \rho^\prime_\varepsilon}}G_{\rho^\prime
_\varepsilon
}\bigr)
\\
&\le&\sup_{t\ge0} e^{-r\Gamma^{\bar{y}}_{t}}\bigl|g\bigl(e^{r\Gamma^{\bar{y}}_{t}}G_{t}
\bigr)\bigr| \le \sup_{t\ge0}e^{-r{t}}\bigl|g\bigl(e^{rt}
\tilde{X}^{\bar{y}}_{t}\bigr)\bigr|.
\end{eqnarray*}
Using this in place of the upper bound on the right-hand side of (\ref
{lebesgue esti}),
we obtain that
\begin{eqnarray*}
\bigl|v\bigl(x,y^\prime\bigr)-v(x,y)\bigr| &\le&3\varepsilon+ \tilde{E}
\bigl(1-e^{-r
(\Gamma_{ A_N^\prime}-N)} \bigr) \sup_{t\ge0}e^{-r{t}}\bigl|g
\bigl(e^{rt}\tilde{X}^{\bar{y}}_{t}\bigr)\bigr|
\\
&&{}+ \tilde{E} e^{-r\Gamma_{ \rho^\prime_\varepsilon}}\bigl |g\bigl(e^
{r\Gamma^\prime_{ \rho^\prime_\varepsilon}}G_{\rho^\prime
_\varepsilon}\bigr)
-g\bigl(e^{r\Gamma_{ \rho^\prime_\varepsilon}}G_{\rho^\prime_\varepsilon}\bigr)\bigr| I\bigl(\rho^\prime_\varepsilon
\le A^\prime_N\bigr).
\end{eqnarray*}
So, after $y$ and $y^\prime$ were replaced by $y_0$ and $y_n$
respectively, it only remains to show that
{
\renewcommand{\theequation}{A.4$'$}
\begin{equation}\label{eq7.4prime}
\lim_{n\to\infty} \tilde{E} e^{-r\Gamma^0_{ \rho^n_\varepsilon}} \bigl|g\bigl(e^{r\Gamma^n_{ \rho^n_\varepsilon}}G_{\rho^n_\varepsilon}
\bigr) -g\bigl(e^{r\Gamma^0_{ \rho^n_\varepsilon}}G_{\rho^n_\varepsilon}\bigr)\bigr| I\bigl(
\rho^n_\varepsilon\le A^n_N\bigr) =
0.%\eqno(\ref{extra term}')
\end{equation}}
This limit refers to the new term in (\ref{extra term})
which was not considered in the proof of Proposition \ref{new int cond}.
But, by dominated convergence, (\ref{eq7.4prime}) would follow
if, for almost every $\omega\in\Omega$, the equality
%
%
%e7.5 #&#
\setcounter{equation}{4}
\begin{eqnarray}
\label{remaining conti}&& \lim_{n\to\infty} \bigl|g\bigl(e^{r\Gamma^n_{ \rho^n_\varepsilon(\omega)}(\omega)}
G_{\rho^n_\varepsilon(\omega)}(\omega)\bigr) -g\bigl(e^{r\Gamma^0_{ \rho^n_\varepsilon(\omega)}(\omega)} G_{\rho^n_\varepsilon(\omega)}(
\omega)\bigr)\bigr|
\nonumber
\\[-8pt]
\\[-8pt]
\nonumber
&&\hspace*{160pt}{}\times I\bigl(\rho^n_\varepsilon(\omega)\le
A^n_N(\omega)\bigr) = 0
\end{eqnarray}
holds, and this is true.
Indeed, choose $\omega\in\Omega$ such that both
$\Gamma^0_{ A_N^n(\omega)}(\omega)\to N$ as $n\to\infty$ and
$t\mapsto G_t(\omega)$ is continuous. Define
\[
c_1 = \sup_{t\le A_N^{\bar{y}_0}(\omega)}\bigl|G_t(\omega)\bigr|,\qquad
c_2 = \Gamma^0_{A_N^{\bar{y}_0}(\omega)}(\omega),
\]
and observe that
\[
0 \le\rho^n_\varepsilon(\omega)I\bigl(\rho^n_\varepsilon(
\omega)\le A^n_N(\omega)\bigr) \le A_N^{\bar{y}_0}(
\omega)I\bigl(\rho^n_\varepsilon(\omega)\le A^n_N(
\omega)\bigr),\qquad n=1,2,\ldots,
\]
since $y_n\downarrow y_0$ and $y_1<\bar{y}_0$ by assumption.
The functions $g$ and $t\mapsto e^{rt}$ are uniformly
continuous on $[-e^{rc_2}c_1,e^{rc_2}c_1]$ and $[0,c_2]$, respectively.
Hence, for the chosen $\omega$, equality (\ref{remaining conti})
follows from
\[
0 \le \bigl(\Gamma^0_{ \rho^n_\varepsilon(\omega)}(\omega) -\Gamma^n_{ \rho^n_\varepsilon(\omega)}(
\omega) \bigr) I\bigl(\rho^n_\varepsilon(\omega)\le
A^n_N(\omega)\bigr) \le \bigl(\Gamma^0_{ A_N^n(\omega)}(
\omega)-N \bigr) \to0
\]
as $n\to\infty$
and almost all $\omega$ are indeed
of this type since the map $t\mapsto G_t$ is
almost surely continuous and
$\lim_{n\to\infty}\Gamma^0_{ A_N^n}$ is almost surely equal to $N$
by Lemma \ref{with feller}.

Part (iii) can be shown
by combining the ideas of the proof of part (ii) and
the proof of Proposition \ref{conti finite time}.
In addition to (\ref{eins}), (\ref{zwei}), (\ref{drei})
there will be an extra term like (\ref{extra term}).
We only need to justify why Lebesgue's dominated convergence theorem
can be applied with respect to this extra term after substituting the
sequence $y_n, n=1,2,\ldots,$ and here,
but only in the case of $y_n\uparrow y_0$,
one needs $g$ to be bounded from below.
\end{pf*}
\end{appendix}

\section*{Acknowledgments}
The authors thank two anonymous referees for their valuable comments.

%
% imsref loaded by akundreckaite, 2014-01-07 08:56:34
%

%


%suskaldyti doi

\printaddresses

\end{document}